\documentclass[11pt]{article}
\usepackage{amssymb}
\usepackage{amsmath}
\usepackage{mathrsfs}
\usepackage{graphics}
\usepackage{graphicx}
\usepackage{subfigure}
\usepackage[T1]{fontenc}
\usepackage{latexsym,amssymb,amsmath,amsfonts,amsthm}\usepackage{txfonts}
\newcommand{\be}{\begin{eqnarray}}
\newcommand{\ben}{\begin{eqnarray*}}
\newcommand{\en}{\end{eqnarray}}
\newcommand{\enn}{\end{eqnarray*}}

\newtheorem{theorem}{Theorem}[section]
\newtheorem{lemma}{Lemma}[section]
\newtheorem{prp}[theorem]{Proposition}
\newtheorem{thm}[theorem]{Theorem}
\newtheorem{cor}[theorem]{Corollary}
\newtheorem{dfn}[theorem]{Definition}

\newtheorem{remark}{Remark}

\begin{document}
\renewcommand{\theequation}{\arabic{section}.\arabic{equation}}
\begin{titlepage}
\title{\bf  Long-time Behavior of 3D
Stochastic Planetary Geostrophic Viscous Model
}
\author{Zhao Dong$^{1,}$,\ \ Rangrang Zhang$^{1,}$\\
{\small $^1$ Institute of Applied Mathematics,
Academy of Mathematics and Systems Science, Chinese Academy of Sciences,}\\
{\small No. 55 Zhongguancun East Road, Haidian District, Beijing, 100190, P. R. China}\\
({\sf dzhao@amt.ac.cn},\ {\sf rrzhang@amss.ac.cn} )}
\date{}
\end{titlepage}
\maketitle

\begin{abstract}
This paper reports the 3D planetary geostrophic viscous model has the exponential ergodicity and global attractor if this model is driven by an additive random noise, which results in the support of the integration of invariant measure for the dynamic is exactly a minimal global attractor. It's worth mentioning that this global attractor has finite Hausdorff dimension.

\end{abstract}
\section{Introduction}
Boussinesq equations is a fundamental model in meteorology, which is deduced from the Primitive Equations using Boussinesq approximation and Hydrostatic approximation (cf. \cite{L-T-W-1,L-T-W-2,L-T-W-3,L-T-W-4}). These are roughly speaking, the Navier-Stokes equations with rotation, forced by the heat buoyancy and coupled with the heat transport equation.
Refer to \cite{STW}, the 3D planetary geostrophic equations are derived from Boussinesq equations using standard scale analysis. In order to strength their regularity, two methods are considered for the momentum equation, one is to add linear drag and the other is to add viscosity term. The latter case is called 3D planetary geostrophic viscous model, which is of our great interest. Compared to 3D Navier-Stokes equations, the important feature of 3D planetary geostrophic viscous model is absence of the nonlinear advection term of the momentum equation, which gives the possibility to obtain its global well-posedness of both weak and strong solution and some other interesting properties.

For the determine case, the early work on this model was derived by Samelson, Termam and Wang in 1998, the
authors established global existence (without uniqueness) of the weak solutions and short time existence of strong solutions in \cite{STW}. In 2003, Cao and Titi  proved the global well posedness and regularity of both the weak and strong solutions to this model in \cite{C-T-2}. Moreover, they show that this model possesses a finite-dimension global attractor. For the random case, in 2016, You and  Li established the existence of global attractor in the space $H\subset L^2$ and proved that the global attractor of the 3D planetary geostrophic viscous model with small linear multiplicative noise converges to the global attractor of the unperturbed case in \cite{YF}.

In this paper, we consider this model in random case from two different aspects: random dynamical system and Markovian semigroup. We obtain its global random attractor in $V\subset H^1$, which is stronger than \cite{YF}. Furthermore, we show that the Hausdorff dimension of the global attractor is almost surely finite. Also we proved that the Markovian semigroup has the exponential mixing property. Thus it can be derived that the support of the invariant measure for the dynamic is a minimal global attractor by the known results.

The main difficulty comes from the nonlinear term $( \int^{-1}_{z}\nabla \cdot vdz')T_{z}$, which results in much more complicated computation. The reason why we can use the method in \cite{O-C} dealing with 2D Navier-Stokes is that the regularity of the velocity can be held by the temperature. This property is obtained by some technical elliptic regularity for nonlocal, stokes-type, second-order elliptic systems in domains with corners. These elliptic regularity results have been used to study the primitive equations in \cite{H-T-Z}.

Since for general multiplicative noise, there is no cocycle property for (\ref{eq-8})-(\ref{eq-14}), which restricts us to dealing with additive noise of the form in Sect. 3.3 (for the linear multiplicative noise, our results also hold). Our results are as follows.

 \begin{thm}\label{thm-1}[\textbf{Global well-posedness}]
 Assume $t_0\in \mathbb{R}$, $\tau\in L^2(D)$, if $T_0\in H(V)$, then for any given $\Upsilon>t_0$, there exists a unique weak (strong) solution $(p_s,v,T)$ of the system (\ref{eq-8})-(\ref{eq-14}) on the interval $[t_0,\Upsilon]$ in the sense of Definition \ref{dfn}; Moveover, the weak (strong) solution $(p_s,v,T)$ is dependent continuously on the initial data.
 \end{thm}
 The definition of space $H$, $V$ and the weak (strong) solution are given in Sect. 3.
 \begin{thm}\label{thm-2}[\textbf{Existence of global attractor}]
   Assume $\tau\in L^2(D)$, then the system (\ref{eq-8})-(\ref{eq-14}) possesses a global random pullback attractor $\mathcal{A}(\omega)$ in $V$. Moreover, $\mathcal{A}(\omega)$ is connected.
 \end{thm}
 The definition of global random attractor is given in Sect.5. Furthermore, we consider the Hausdorff dimension of $\mathcal{A}(\omega)$.
 \begin{thm}\label{thm-8}[\textbf{Finite Hausdorff dimension of the global attractor}]
 Assume $\tau\in L^2(D)$, then the Hausdorff dimension of the global attractor $\mathcal{A}(\omega)$ is finite, $\mathbb{P}-$ a.s.
 \end{thm}
 Ergodic property of the 3D stochastic planetary geostrophic viscous model is also established. At that time, some nondegenerate condition on the noise is needed.
 \begin{thm}\label{thm-3}[\textbf{Exponential mixing}]
   Assume $\tau\in L^2(D)$ and \textbf{Hypothesis H0} holds, then there exists a unique stationary probability measure $\mu$ of $(\mathcal{P}_t)_{t\in \mathbb{R}^{+}}$ on $H$. Moreover, $\mu$ satisfies
 \begin{equation}\notag
 \int_{H}|x|^2\mu(dx)<\infty,
 \end{equation}
 and there exist C, $\gamma^{'}>0 $ such that for any $\lambda\in\mathcal{P}(H)$
 \begin{equation}\notag
 \|\mathcal{P}^*_t\lambda-\mu\|_{*}\leq Ce^{-\gamma^{'}t}\left(1+\int_{H}|x|^2\lambda(dx)\right).
 \end{equation}
 \end{thm}
where $\mathcal{P}(H)$ stands for the set of all probability measures on $H$. The notation of the norm $\|\cdot\|_{*}$ and \textbf{Hypothesis H0} are given in Sect. 6.

This paper is organized as follows: In Sects. 2 and 3, we introduce the 3D stochastic planetary geostrophic viscous model, and review some basic representational results; The global well-posedness of both weak and strong solution are proved in Sect. 4; In Sect. 5, the existence of global attractor for the strong solution is obtained; In Sect. 6, the exponential mixing for the weak solution is established; Finally, in Sect. 7, we discuss connection between global attractor and the invariant measure of this model.

\section{Preliminaries}
The 3D Planetary Geostrophic Viscous Model under a stochastic forcing, in a Cartesian system, are written as
\begin{eqnarray}\label{eq-1}
\nabla P+f{k}\times v + \varepsilon L_1 v=0,\\
\label{eq-2}
\partial_{z}P+T=0,\\
\label{eq-3}
\nabla\cdot v+\partial_{z}\theta=0,\\
\label{eq-4}
\frac{\partial T}{\partial t}+(v\cdot\nabla)T+\theta\frac{\partial T}{\partial z}+L_2 T=\Theta(t,x,y,z),
\end{eqnarray}
where the horizontal velocity field $v=(v_{1},v_{2})$, the three-dimensional velocity field $(v_{1},v_{2},\theta)$, the temperature $T$ and the pressure $P$ are unknown functionals. $f$ is the Coriolis parameter and ${k}$ is vertical unit vector. $\varepsilon L_1 v$ is the viscous term, where $\varepsilon$ is a positive viscous constant. $\Theta(t,x,y,z)$ is the stochastic forcing which will be described in Sect.3.

Set $\nabla=(\partial x,\partial y)$ to be the horizontal gradient operator and  $\Delta=\partial^{2}_{x}+\partial^{2}_{y}$  to be the horizontal Laplacian. The viscosity and the heat diffusion operators $L_1$ and $L_2$ are given by
\begin{eqnarray*}
L_1v&=&-A_h\Delta v -A_v\frac{\partial^2 v}{\partial z^2},\\
L_2T&=&-K_h\Delta T -K_v\frac{\partial^2 T}{\partial z^2},
\end{eqnarray*}
where $A_h$, $A_v$ are positive molecular viscosities and $K_h$, $K_v$ are positive conductivity constants.

The space domain for (\ref{eq-1})-(\ref{eq-4}) is
\[
\mathcal{O}=\textit{D}\times(-1,0),
\]
where $\textit{D}$ is a  bounded open domain with smooth bandary in $\mathbb{R}^{2}$.

The boundary conditions for (\ref{eq-1})-(\ref{eq-4}) are given by
 \begin{eqnarray}
 \label{eq-5}
 A_v\partial_{z}v=\tau,\ \theta=0,\ K_v\partial_{z}T=-\alpha(T-T^*) & {\rm on}& \textit{D}\times\{0\}=\Gamma_{\textit{u}},\\
 \label{eq-6}
 \partial_{z}v=0,\ \theta=0,\ \partial_{z}T=0     &{\rm on}&  \textit{D}\times\{-1\}=\Gamma_{\textit{b}},\\
 \label{eq-7}
 v\cdot {n}=0,\  \frac{\partial v}{\partial n}\times {n}=0,\ \frac{\partial T}{\partial {n}}=0                                    &{\rm on}& \partial \textit{D}\times (-1,0)=\Gamma_{l},
 \end{eqnarray}
 where $\tau(x,y)$ is the wind stress, ${n}$ is the normal vector to $\Gamma_{l}$, $T^*(x,y)$ is the typical temperature of the top surface, and $\alpha$ is a positive constant.

Refer to \cite{C-T-2}, $\tau$ and $T^*$ are assumed to satisfy the compatibility boundary conditions
\begin{eqnarray}\label{equ-2}
\tau\cdot n=0,\quad \frac{\partial \tau}{\partial n}\times n=0, \quad \rm{on}\ \partial D.
\end{eqnarray}
\begin{eqnarray}\label{equ-3}
\frac{\partial T^{*}}{\partial n}=0 \quad \rm{on}\ \partial D.
\end{eqnarray}
In addition, the initial condition is
\begin{eqnarray}\label{equ-1}
T(x,y,z,t_0)=T_0(x,y,z).
\end{eqnarray}
\begin{remark}
For simplicity and without loss generality we will assume $\varepsilon,A_h,A_h,K_h,K_v$ are equal to 1. Moreover, By \cite{C-T-1}, we know that due to the compatibility boundary conditions (\ref{equ-2}) and (\ref{equ-3}), one  can convert the boundary condition (\ref{eq-5})-(\ref{eq-7}) to be homogeneous by replacing $(v,T)$
by $(v+\frac{(z+1)^2-1/3}{2}\tau,T+T^*)$. For convenience, in this article, we suppose $ T^*=0$ and we emphasize that our results are still valid for general $T^*$ provided it is smooth enough.
\end{remark}
Let us reformulate the system (\ref{eq-1})-(\ref{equ-1}).
By integrating (\ref{eq-3}) from -1 to $\textit{z}$ and using (\ref{eq-5}), (\ref{eq-6}), we obtain
\begin{equation}\label{eq-90}
\theta(t,x,y,z)=\Phi(v)(t,x,y,z)=-\int^{z}_{-1}\nabla\cdot v(t,x,y,z')dz'.
\end{equation}
Integrating (\ref{eq-2}) from $-1$ to $z$, we obtain
\begin{equation}\label{eq-91}
P(x,y,z,t)=-\int^z_{-1}T(x,y,z',t)dz'+p_s(x,y,t),
\end{equation}
where $p_{s}(x,y,t)$ is a certain unknown function at $\Gamma_{b}$.

Under the above calculation, (\ref{eq-1})-(\ref{eq-4}) can be rewritten as
\begin{eqnarray}\label{eq-8}
\nabla p_{s}-\int^{z}_{-1}\nabla T dz'+f{k}\times v + L_1v &=& 0,\\
\label{eq-9}
\frac{\partial T}{\partial t}+(v\cdot\nabla)T+\Phi(v)\frac{\partial T}{\partial z}+ L_2T&=&\Theta(t,x,y,z),\\
\label{eq-10}
\int^{0}_{-1}\nabla\cdot v dz&=&0.
\end{eqnarray}
The boundary  and initial conditions are
\begin{eqnarray}\label{eq-11}
 \partial_{z}v=\tau,\ \partial_{z}T=-\alpha T     &{\rm on}& \Gamma_{\textit{u}},\\
 \label{eq-12}
 \partial_{z}v=0,\ \partial_{z}T=0      &{\rm on}& \Gamma_{\textit{b}},\\
 \label{eq-13}
 v\cdot n=0,\ \frac{\partial v}{\partial n}\times n=0,\ \frac{\partial T}{\partial {n}}=0      &{\rm on}& \Gamma_{l},\\
\label{eq-14}
T(x,y,z,t_0)=T_0(x,y,z).
\end{eqnarray}
\section{Functional Spaces and Inequalities }
\subsection{Some Functional Spaces }
Let $\mathcal{L}(K_1,K_2)$ (resp. $\mathcal{L}_2(K_1,K_2)$) be the space of bounded (resp. Hilbert-Schmidt) linear operators from the Hilbert space $K_1$ to $K_2$, whose norms are denoted by $\|\cdot\|_{\mathcal{L}(K_1,K_2)}$ ( resp $\|\cdot\|_{\mathcal{L}_2(K_1,K_2)}$). Denote by $|\cdot|_{L^2(D)}$ the norm of $L^2(D)$. $|\cdot|$ and $(\cdot,\cdot)$ represent the norm and inner product in $L^2(\mathcal{O})$. Let $|\cdot|_{p}$ and $|\cdot|_{H^{p}(\mathcal{O})}$ be the norms of $L^p(\mathcal{O})$ and Sobolev space $H^{p}(\mathcal{O})$, respectively, for integer $p$,
\begin{equation}\notag
\left\{
  \begin{array}{ll}
    H^p(\mathcal{O})=\Big\{Y\in L^2(\mathcal{O})\Big| \ \partial_{\alpha}Y\in L^2(\mathcal{O}),\  \ 0\leq|\alpha|\leq p\Big\},&  \\
    |Y|^2_{H^p(\mathcal{O})}=\sum_{0\leq|\alpha|\leq p}|\partial_{\alpha}Y|^2. &
  \end{array}
\right.
\end{equation}

Now, define working spaces for (\ref{eq-8})-(\ref{eq-14}). Let
 \begin{eqnarray}\notag
 &&\breve{\mathcal{V}}\triangleq\left\{v\in (C^{\infty}(\mathcal{O}))^2;\  \frac{\partial v}{\partial z}\Big|_{\Gamma_u}=\tau,\ \frac{\partial v}{\partial z}\Big|_{\Gamma_b}=0,\  v\cdot{n}\Big|_{\Gamma_l}=0,\  \frac{\partial v}{\partial {n}}\times {n}\Big|_{\Gamma_l}=0,\ \int^{0}_{-1}\nabla \cdot v dz=0\right\},\\ \notag
 &&{\mathcal{V}}\triangleq\left\{T\in C^{\infty}(\mathcal{O}); \frac{\partial T}{\partial z}+\alpha T\Big|_{\Gamma_u}=0, \frac{\partial T}{\partial z}\Big|_{\Gamma_b}=0,  \frac{\partial T}{\partial {n}}\Big|_{\Gamma_l}=0 \right\},
 \end{eqnarray}
$\breve{V}$= the closure of $\breve{\mathcal{V}} $ with respect to the norm $|\cdot|_{H^1(\mathcal{O})}$,\\
$\breve{H}$= the closure of $\breve{\mathcal{V}} $ with respect to the norm $|\cdot|$,\\
$V$= the closure of ${\mathcal{V}} $ with respect to the norm $\|\cdot\|$,\\
$H$= the closure of ${\mathcal{V}} $ with respect to the norm $|\cdot|$,

where
\begin{eqnarray*}
&&|v|^2_{H^1(\mathcal{O})}=\int_{\mathcal{O}}|\nabla v|^2dxdydz+\int_{\mathcal{O}}|\frac{\partial v}{\partial z}|^2dxdydz,\quad\quad\quad\quad\quad\quad\quad\quad\quad\quad\quad\quad\quad\quad\quad\quad\\
&&|v|^2=\int_{\mathcal{O}}|v|^2dxdydz,\\
&&\|T\|^2=\int_{\mathcal{O}}|\nabla T|^2dxdydz+\int_{\mathcal{O}}|\frac{\partial T}{\partial z}|^2dxdydz+\alpha\int_{\Gamma_u}|T(x,y,0)|^2dxdy,\\
&&|T|^2=\int_{\mathcal{O}}|T|^2dxdydz.
\end{eqnarray*}
Define the bilinear functional $a:V\times V\rightarrow \mathbb{R}$, and its corresponding linear operator $A: V\rightarrow V^{'}$, by
\begin{eqnarray*}
a(T,T_1)\triangleq(AT, T_1)=\int_{\mathcal{O}}\left(\nabla T\cdot \nabla T_1+\frac{\partial T}{\partial z}\frac{\partial T_1}{\partial z}\right)dxdydz+\alpha\int_{\Gamma_u}TT_1dxdy,
\end{eqnarray*}
for any $T,T_1\in V$.

It is well known that the operator $A$ with domain
\[
D(A)=closure\ of\ \mathcal{V}\ with\ respect\ to\ the \ H^2(\mathcal{O})\ topology.
\]
is a positive self-adjoint operator with discrete spectrum in $H$. Let $(e_n)_{n=1,2,\cdot\cdot\cdot}$ be an eigenbasis of $H$ associated to the increasing sequence $\{\lambda_n\}_{n=1,2,\cdot\cdot\cdot}$ of eigenvalues of $A$. Denote by $H_n=span\{e_1,\cdot\cdot\cdot,e_n\}$, and $P_n: H\rightarrow H_n$ the $H$ orthogonal projection onto $H_n$. Moreover, by definition of $\breve{H}$, there exists an orthogonal basis $(r_n)_{n=1,2,\cdot\cdot\cdot}$ of $\breve{H}$.

\begin{lemma}(\cite{C-T-2}) We have the following results:
\begin{description}
  \item[(i)] There exists a positive constant $K_1$ such that for $T\in H^1(\mathcal{O})$,
\[
\frac{1}{K_1}\|T\|^2\leq |T|^2_{H^1(\mathcal{O})}\leq K_1\|T\|^2.
\]
\item[(ii)] For every $T\in V$, $\|T\|^2=(AT,T)$.
\item[(iii)] For every $T\in V$,
  $|T|^2_{H^2(\mathcal{O})}\leq (AT,AT)$.

\end{description}
\end{lemma}

\subsection{Some Inequalities}
We firstly recall some interpolation inequalities which will be used later (see \cite{Adams} for detail).\\
For $h\in H^1(D)$,
\begin{eqnarray*}
|h|_{L^4(D)}&\leq& c|h|^{\frac{1}{2}}_{L^2(D)}|h|^{\frac{1}{2}}_{H^1(D)},\\
|h|_{L^5(D)}&\leq& c|h|^{\frac{3}{5}}_{L^3(D)}|h|^{\frac{2}{5}}_{H^1(D)},\\
|h|_{L^6(D)}&\leq& c|h|^{\frac{2}{3}}_{L^4(D)}|h|^{\frac{1}{3}}_{H^1(D)}.
\end{eqnarray*}
For $h\in H^1(\mathcal{O})$,
\begin{eqnarray*}
|h|_3&\leq& c|h|^{\frac{1}{2}}|h|^{\frac{1}{2}}_{H^1(\mathcal{O})},\\
|h|_6&\leq& c|h|_{H^1(\mathcal{O})},\\
|h|_{\infty}&\leq& c|h|^{\frac{1}{2}}_{H^1(\mathcal{O})}|h|^{\frac{1}{2}}_{H^2(\mathcal{O})}.
\end{eqnarray*}

\begin{lemma}(\cite{C-T-2})\label{lem-5}
Suppose $u=(u_1,u_2)$ be a smooth vector field, and let $f$ and $g$ be smooth scalar functions, then
\begin{eqnarray}\label{eq-34}
\Big|\int_{\mathcal{O}}(u\cdot \nabla) f gdxdydz\Big|&\leq &|u|_6|\nabla f|_3|g|\leq c|u|_{H^1(\mathcal{O})}|\nabla f|^{\frac{1}{2}}|f|^{\frac{1}{2}}_{H^2(\mathcal{O})}|g|,\\
\label{eq-37}
\Big|\int_{\mathcal{O}}(\nabla \cdot \int^{z}_{-1}u(x,y,z',t)dz')fgdxdydz \Big|&\leq & c|f||u|^{\frac{1}{2}}_{H^1(\mathcal{O})}|u|^{\frac{1}{2}}_{H^2(\mathcal{O})}\|g\|^{\frac{1}{2}}|g|^{\frac{1}{2}}.
\end{eqnarray}
\end{lemma}

\

Consider the following equations
\begin{eqnarray}\label{equ-27}
\nabla p_{s}(x,y,t)-\int^{z}_{-1}\nabla g dz'+f{k}\times \xi + L_1\xi &=& 0,\\
\label{equ-28}
\int^{0}_{-1}\nabla\cdot \xi  dz&=&0,\\
\label{equ-29}
\frac{\partial \xi}{\partial z}\Big|_{z=0}=0,\  \frac{\partial \xi}{\partial z}\Big|_{z=-1}=0, \  \xi\cdot n\Big|_{\Gamma_{l}}=0,\ \frac{\partial \xi}{\partial n}\times n\Big|_{\Gamma_{l}}&=&0.
\end{eqnarray}

The following  proposition  plays an important role during the following sections.

\begin{prp}\label{prp-2}
For  (\ref{equ-27})-(\ref{equ-29}), we have
\begin{equation}\label{eq-26}
|\xi|_{H^1(\mathcal{O})}\leq c |g|,\quad |\xi|_{H^2(\mathcal{O})}\leq c \|g\|.
\end{equation}
\end{prp}
\textbf{Proof of Proposition \ref{prp-2}} \quad By averaging (\ref{equ-27}) and (\ref{equ-29}) from $-1$ to 0 with respect to $z$, we obtain
\begin{eqnarray}\label{eq-19-1}
\nabla p_s-\int^0_{-1}(\int^z_{-1}\nabla g(x,y,z',t)dz')dz+fk\times\bar{\xi}-\Delta\bar{\xi}-\int^0_{-1}\frac{\partial^2 \xi}{\partial z^2}dz=0,\\
\label{eq-20}
\nabla\cdot \bar{\xi}=0,\\
\label{eq-21}
\bar{\xi}\cdot n=0, \quad \frac{\partial \bar{\xi}}{\partial n}\times n =0 \quad on\ \partial D,
\end{eqnarray}
where
\[
\bar{\xi}(x,y,t)=\int^{0}_{-1}\xi(x,y,z,t)dz.
\]
From the Fubini Theorem and boundary conditions (\ref{equ-29}), (\ref{eq-19-1}) can be written as
\begin{eqnarray}\label{eq-19}
f{k}\times \bar{\xi}+\nabla p_s+\nabla\left(\int^{0}_{-1}z'g(x,y,z',t)dz'\right)-\Delta\bar{\xi}=0.
\end{eqnarray}
By taking inner product on both side  of (\ref{eq-19}) with\ $\bar{\xi}$ in $L^{2}{(\mathcal{O})}$, we obtain
 \[
 \int_{\mathcal{O}}\left[\nabla\left(p_s(x,y,t)+\int^{0}_{-1}z'g(x,y,z',t)dz'\right)-\Delta \bar{\xi}\right]\bar{\xi}dxdydz=0.
 \]
Notice that (\ref{equ-28}) and (\ref{equ-29}), by taking integration by parts, we get
 \[
 \int_{\mathcal{O}}|\nabla \bar{\xi}|^2dxdydz=0.
 \]
 Thus, $\bar{\xi}$ is a constant function. By (\ref{equ-29}), we reach $\bar{\xi}=0$. As a result, we have
 \[
 p_s(x,y,t)=-\int^{0}_{-1}z'g(x,y,z',t)dz',
 \]
 ($p_s$ is unique up to a constant), then, (\ref{equ-27}) can be written as
 \begin{equation}\notag
-\nabla [\int^{0}_{-1}z'g(x,y,z',t)dz'+\int^{z}_{-1}g(x,y,z',t)dz']+f{k}\times \xi+L_1 \xi=0.
 \end{equation}
 From the boundary value problem of the second-order elliptic equation, we have the following regularity results (refer to similar techniques to those developed in \cite{HTZ},\cite{Z}):
\begin{eqnarray}\notag
|\xi|_{H^1(\mathcal{O})}\leq c |g|,\quad |\xi|_{H^2(\mathcal{O})}\leq c \|g\|.
\end{eqnarray}$\hfill\blacksquare$
\subsection{Stochastic Forcing}
In the present paper, we assume that the stochastic forcing $\Theta(t,x,y,z)$ is an additive white noise with the form
\[
\Theta(t,x,y,z)=G\frac{dW(t)}{dt},
\]
where $G$ is a Hilbert-Schmidt operator from $H$ to $H^{1}(\mathcal{O})$, i.e.
\begin{eqnarray}\label{equ-8}
\sum^{+\infty}_{i=1}\|Ge_i\|^2=\sum^{+\infty}_{i=1}\lambda_i|Ge_i|^2<\infty,
\end{eqnarray}
and the random process $W$ is a two-sided in time cylindrical wiener process in $H$ with the form $W(t)=\sum^{+\infty}_{i=1}w_i(t,\omega)e_i$. Here, $\{w_i\}_{i\geq 1}$ is a sequence of independent standard one-dimensional real-valued Brownian motions on a complete probability space $(\Omega, \mathcal{F}, \mathbb{P})$ with expectation $\mathbb{E}$. More precisely, let
\[
\Omega=\Big\{\omega: \omega\in C(\mathbb{R},\mathbb{R}^{\infty}), \omega(0)=0\Big\},
\]
with $\mathbb{P}$ being a product measure of two Wiener measures on the negative and positive time parts of $\Omega$.
\subsection{An Auxiliary Ornstein-Uhlenbeck Process}
For $\alpha\geq 0$, let
\[
Z_{\alpha}(t)=\int^t_{-\infty}e^{(t-s)(-A-\alpha)}GdW(s),
\]
be the solution of the stochastic Stokes equation
\[
dZ=(-\alpha Z-AZ)dt+GdW(t),
\]
with $Z(0)=\int^0_{-\infty}e^{s(A+\alpha)}GdW(s)$. In the following, denote $Z_{\alpha}(t)$ by $Z(t)$ for simplicity. The damping term $\alpha Z$ is not necessary to the global well-posedness, but is very useful to prove the existence of global attractor.
\begin{lemma}\label{Guo}\cite{Guo}
If $Z(t)$ is the solution for the above equation, then the process $Z(t)$ is a stationary ergodic solution with continuous trajectories in $D(A)$.
\end{lemma}
\begin{prp}\label{prp-5}
For any $t$, $\mathbb{E}|AZ(t)|^2<\infty$ and $\mathbb{E}|AZ(t)|^2\rightarrow 0$, as $\alpha\rightarrow +\infty$.
\end{prp}
\textbf{Proof of Proposition \ref{prp-5}} \quad We prove the second statement, the first one can be proved  similarly.
\begin{eqnarray}\label{equ-9}
\mathbb{E}|AZ(t)|^2&=&\mathbb{E}\sum^{+\infty}_{i=1}|\int^t_{-\infty}\lambda_ie^{(t-s)(-\lambda_i-\alpha)}Ge_idw_i(s)|^2\\ \notag
&=&\sum^{+\infty}_{i=1}\int^t_{-\infty}\lambda^2_ie^{2(t-s)(-\lambda_i-\alpha)}|Ge_i|^2ds\\ \notag
&=&\sum^{+\infty}_{i=1}\frac{\lambda^2_i|Ge_i|^2}{2(\lambda_i+\alpha)}, \notag
\end{eqnarray}
by (\ref{equ-8}), we obtain the result. $\hfill\blacksquare$
\begin{remark}
A simple example for $\Theta(t,x,y,z)$ is
\[
\Theta(t,x,y,z)=\sum^{N}_{j=1}\mu_j\frac{dw_j(t)}{dt} \quad  for \ some \  N.
\]
 More examples can be found in \cite{Guo}.
\end{remark}

\subsection{Definition of Weak (Strong) Solution of (\ref{eq-8})-(\ref{eq-14})}
\begin{dfn}\label{dfn}
 Assume $t_0\in \mathbb{R}$ and ${\Upsilon}>t_0$, a stochastic process $(p_s, v, T)$ is called a weak solution of (\ref{eq-8})-(\ref{eq-14}) in $[t_0,\Upsilon]$, if \ $\mathbb{P}-$ a.e. $\omega\in \Omega$
 \begin{eqnarray*}
p_s(x,y,t)&\in& C(t_0, \Upsilon; L^2(D)) \bigcap L^2(t_0,\Upsilon; H^1(D)),\\
v(x,y,z,t)&\in& C(t_0,\Upsilon; H^1(\mathcal{O}))\bigcap L^2(t_0,\Upsilon; H^2(\mathcal{O})),\\
T(x,y,z,t)&\in& C(t_0,\Upsilon; H)\bigcap L^2(t_0,\Upsilon; V),
\end{eqnarray*}
and for all $(\phi, \psi) \in H^2(\mathcal{O})\times D(A)$, the identities hold $\mathbb{P}-a.s.$
\begin{eqnarray*}
-( p_s,\nabla \cdot \phi)+\Big(\int^{z}_{-1} T dz',\nabla \cdot \phi\Big)+(fk\times v, \phi)+(\nabla v,\nabla \phi)+\Big(\frac{\partial v}{\partial z},\frac{\partial \phi}{\partial z}\Big)&=&\int_{\Gamma_u}\tau\phi dxdydz,\\
\Big(T(t),\psi\Big)-\int^{t}_{t_0}\Big((v\cdot \nabla) \psi, T\Big)ds - \int^{t}_{t_0}\Big(\Phi(v)\frac{\partial \psi}{\partial z}, T\Big)ds +\int^{t}_{t_0}(A\psi, T)ds
&=&(T_0, \psi)+ \Big(\int^{t}_{t_0}GdW(s),\psi\Big).
\end{eqnarray*}

Moreover, a weak solution is called a strong solution of (\ref{eq-8})-(\ref{eq-14}) on $[t_0,\Upsilon]$, if in addition, it satisfies
\begin{eqnarray*}
p_s(x,y,t)&\in& C(t_0, \Upsilon; H^1(D)) \bigcap L^2(t_0,\Upsilon; H^2(D)),\\
v(x,y,z,t)&\in& C(t_0,\Upsilon; H^1(\mathcal{O}))\bigcap L^2(t_0,\Upsilon; H^2(\mathcal{O})),\\
T(x,y,z,t)&\in& C(t_0,\Upsilon; V)\bigcap L^2(t_0,\Upsilon; D(A)).
\end{eqnarray*}
\end{dfn}

 \section{Global Well-posedness}
We firstly prove the global well-posedness of weak solution.

\begin{thm}\label{thm-4}(Global well-posedness of weak solution)
Suppose that $\tau\in L^2(D)$. Then for every $T_0\in H$, there is a unique weak solution $(p_s,v,T)$ of the system (\ref{eq-8})-(\ref{eq-14}) in the sense of definition \ref{dfn}.
\end{thm}
\textbf{Proof of Theorem \ref{thm-4} } \quad \textbf{(Existence of weak solution)} \quad Let $h(t)=T(t)-Z(t)$, where $Z(t)$ is defined in Sect.3.4, we get the following random parameter equations:
\begin{eqnarray}\label{eq-35}
\frac{dh}{dt}+L_2h+(v\cdot\nabla)(h+Z)+\Phi(v)\frac{\partial(h+Z)}{\partial z}&=&\alpha Z.\\
\label{equ-10}
\nabla\Big[p_s-\int^{z}_{-1} (h+Z) dz'\Big]+fk\times v +L_1 v&=&0,\\
\label{equ-11}
\frac{\partial v}{\partial z}\Big|_{z=0}=\tau, \quad  \frac{\partial v}{\partial z}\Big|_{z=-1}=0, \quad v\cdot n\Big|_{\Gamma_l}=0,\quad \frac{\partial v}{\partial n}\times n\Big|_{\Gamma_l}&=&0,\\
\label{equ-12}
h(0)&=&T_0-Z_{t_0}.
\end{eqnarray}
In the above system, the unknowns are $h(x,y,z,t)$, $p_s(x,y,t)$ and $v(x,y,z,t)$, while $T_0$ and $Z_{t_0}$ are given. It's clear that once we determine $h,p_s$ and $v$, we can easily recover the original unknowns $(v,\theta)$, $P$ and $T$ of the system (\ref{eq-1})-(\ref{eq-4}) by (\ref{eq-90}) and (\ref{eq-91}).

In the following, we proceed by using Galerkin method to prove the global existence of weak solution using the similar method as \cite{STW}. We divide the proof into several steps as follows.

\textbf{\emph{Step 1. Approximate Solutions.}}
Let $\omega\in \Omega$ be fixed. In order to obtain a solution of (\ref{eq-35})-(\ref{equ-12}), we find an approximate solution pair $(h_n, v_n)$ of the form
\begin{eqnarray*}
h_n(x,y,z,t)&=&\sum^n_{k=1}\beta^n_k(t)e_k(x,y,z),\\
v_n(x,y,z,t)&=&\sum^n_{k=1}\gamma^n_k(t)r_k(x,y,z),
\end{eqnarray*}
where $\{\beta^n_k\}$ and $\{\gamma^n_k\} $ are real functionals and  $(h_n, v_n)$
satisfies
\begin{eqnarray}\label{eq-24}
\frac{d }{dt}(h_n, e_k)+ \Big((v_n\cdot \nabla)(h_n+Z_n),e_k\Big) + \Big( \Phi(v_n)\frac{\partial (h_n+Z_n)}{\partial z} , e_k\Big)+(L_2 h_n,e_k)&=&(\alpha Z_n,e_k),\\
\label{eq-25}
-\Big(\int^{z}_{-1}\nabla (h_n+Z_n) dz', r_k\Big)+(fk\times v_n,r_k) +(L_1 v_n, r_k)&=&0,\\
\label{eq-25-2}
h_n(0)&=&P_n(T_0-Z_{t_0}),
\end{eqnarray}
for $k=1,\cdot\cdot\cdot,n$, where $Z_n=P_n Z$ and $(\nabla p_s, r_k)=0$ are used.

It's easy to see that (\ref{eq-24}) and (\ref{eq-25-2}) are equivalent to an initial value problem of a system of $n$ ODEs, and (\ref{eq-25}) is equivalent to a linear system of $n$ equations. Hence the existence of the local (in time) approximate solution of (\ref{eq-24})-(\ref{eq-25-2}) is obvious. For a short interval of time $[t_0, \Upsilon^*)$, where $\Upsilon^*$ is independent of $n$, we have the existence and uniqueness of $h_n(x,y,z,t)$ and $v_n(x,y,z,t)$, then $p_s(x,y,t)$ determined by (\ref{equ-10}) is also obtained on $[t_0, \Upsilon^*)$.

\textbf{\emph{Step 2. Energy estimates.}} To obtain the global (in time) solution of the original problem, the energy estimates are necessary.

Firstly, multiplying (\ref{eq-25}) by $\gamma^n_k(t)$ and adding up the resulting equations for $k=1,\cdot\cdot\cdot,n$, we deduce that
\[
(-\int^z_{-1}\nabla (h_n+Z_n)dz',v_n)+(fk\times v_n, v_n)+(L_1v_n, v_n)=0.
\]
Since
\[
(L_1 v_n,v_n)=|v_n|^2_{H^1(\mathcal{O})}-\int_{\Gamma_u}v_n \tau dxdy, \quad (fk\times v_n, v_n)=0,
\]
we have
\begin{eqnarray*}
|v_n|^2_{H^1(\mathcal{O})}&=&\int_{\Gamma_s}v_n \tau dxdy+(-\int^z_{-1}\nabla (h_n+Z_n)dz',v_n)\\
&\leq& c|h_n+Z_n||\nabla v_n|+c|\tau|_{L^2(D)}|v_n|_{L^2(\Gamma_u)}\\
&\leq& \frac{1}{2}|\nabla v_n|^2+c(|h_n|^2+|Z_n|^2+|\tau|_{L^2(D)}^2),
\end{eqnarray*}
thus,
\begin{eqnarray}\label{eq-80}
|v_n|^2_{H^1(\mathcal{O})}\leq c(|h_n|^2+|Z_n|^2+|\tau|_{L^2(D)}^2).
\end{eqnarray}
For $|p_s|^2_{H^1(D)}$, by
\[
(\nabla p_s,L_1 v_n)=-\int_{\Gamma_u}\tau \nabla p_s dxdy
\leq c|\tau|_{L^2(D)}|p_s|_{H^1(D)},
\]
we have
\begin{eqnarray*}
|p_s|^2_{H^1(D)}&=&(\nabla p_s,\nabla p_s)\\
&=&(\nabla p_s, \int^{z}_{-1} (h_n+Z_n) dz'-fk\times v_n -L_1 v_n)\\
&\leq& c\|h_n+Z_n\||p_s|_{H^1(D)}+c|v_n||p_s|_{H^1(D)}+c|\tau|_{L^2(D)}|p_s|_{H^1(D)},
\end{eqnarray*}
therefore,
\begin{eqnarray}\label{eq-82}
|p_s|^2_{H^1(D)}\leq c(\|h_n\|^2+\|Z_n\|^2+|\tau|_{L^2(D)}^2).
\end{eqnarray}
For $|L_1v_n|^2$, we have
\begin{eqnarray*}
|L_1v_n|^2&=&(\int^z_{-1}\nabla (h_n+Z_n)dz',L_1v_n)-(fk\times v_n, L_1 v_n)-(\nabla p_s,L_1 v_n)\\
&\leq& c|\tau|_{L^2(D)}|p_s|_{H^1(D)}+c|L_1 v_n||v_n|+c|L_1 v_n|\|h_n+Z_n\|\\
&\leq & \frac{1}{2}|L_1 v_n|^2+c(\|h_n\|^2+\|Z_n\|^2+|\tau|_{L^2(D)}^2),
\end{eqnarray*}
where (\ref{eq-82}) is used. Then, we get
\begin{eqnarray}\label{eq-81}
|v_n|^2_{H^2(\mathcal{O})}\leq|L_1 v_n|^2\leq c(\|h_n\|^2+\|Z_n\|^2+|\tau|_{L^2(D)}^2).
\end{eqnarray}
The above (\ref{eq-80})-(\ref{eq-82}) hold for any $t$ in the interval of the maximal existence of the solutions of the approximate problem (\ref{eq-24})-(\ref{eq-25-2}).

Multiplying (\ref{eq-24}) by $\beta^n_k(t)$ and adding the resulting equations up, we get
\begin{eqnarray*}
\frac{1}{2}\frac{d |h_n|^2}{dt}+\|h_n\|^2+ \Big([(v_n\cdot \nabla)(h_n+Z_n)+ \Phi(v_n)\frac{\partial (h_n+Z_n)}{\partial z}],h_n\Big)=(\alpha Z_n, h_n),
\end{eqnarray*}
by integration by parts and the boundary conditions (\ref{eq-11})-(\ref{eq-13}), we have
\[
\Big([(v_n\cdot \nabla)h_n+ \Phi(v_n)\frac{\partial h_n}{\partial z}],h_n\Big)=0,
\]
furthermore, by H\"{o}lder inequality we have
\[
|\Big((v_n\cdot \nabla)Z_n,h_n \Big)|
\leq c|\nabla Z_n||h_n|_3|v_n|_6,
\]
by (\ref{eq-80}), we have
\[
|v_n|_6\leq c|v_n|_{H^1(\mathcal{O})}\leq c(|h_n|+|Z_n|+|\tau|_{L^2(D)}),
\]
then, by the Young inequality, we have
\begin{eqnarray*}
|\Big((v_n\cdot \nabla)Z_n,h_n \Big)|
&\leq& c|\nabla Z_n|\|h_n\|(|h_n|+|Z_n|+|\tau|_{L^2(D)})\\
&\leq &\varepsilon \|h_n\|^2+c|\nabla Z_n|^2|h_n|^2+c|\nabla Z_n|^2(|Z_n|^2+|\tau|^2_{L^2(D)}).
\end{eqnarray*}
Moreover, by Lemma \ref{lem-5}, we obtain
\begin{eqnarray*}
|\int_{\mathcal{O}}\Phi(v_n)\frac{\partial Z_n}{\partial z}h_n dxdydz|
&\leq &c|\frac{\partial Z_n}{\partial z}||\nabla v|^{\frac{1}{2}}|v|^{\frac{1}{2}}_{H^2(\mathcal{O})}| h_n|^{\frac{1}{2}}\|h_n\|^{\frac{1}{2}}\\
&\leq& c|\frac{\partial Z_n}{\partial z}|(| h_n|^{\frac{1}{2}}\|h_n\|^{\frac{1}{2}}+\|Z_n\|+|\tau|_{L^2(D)})| h_n|^{\frac{1}{2}}\|h_n\|^{\frac{1}{2}}\\
&\leq& c|\frac{\partial Z_n}{\partial z}|| h_n|\|h_n\|+c|\frac{\partial Z_n}{\partial z}|(\|Z_n\|+|\tau|_{L^2(D)})| h_n|^{\frac{1}{2}}\|h_n\|^{\frac{1}{2}}\\
&\leq& \varepsilon \|h_n\|^2+c(|\frac{\partial Z_n}{\partial z}|^2+|\frac{\partial Z_n}{\partial z}|^4\|Z_n\|^4+|\frac{\partial Z_n}{\partial z}|^4|\tau|^4_{L^2(D)})| h_n|^2.
\end{eqnarray*}
where (\ref{eq-80}) and the Young inequality are used in the second and the forth inequality, respectively.

Applying Cauchy-Schwarz inequality, we obtain
\begin{eqnarray*}
|\int_{\mathcal{O}}\alpha Z_n h_n dxdydz|
&\leq& c|Z_n||h_n|\\
&\leq& c|Z_n|^2+c|h_n|^2.
\end{eqnarray*}
Therefore, collecting all the above inequalities, we reach
\begin{eqnarray}\label{eq-23}
\frac{d |h_n|^2}{dt}+\|h_n\|^2
&\leq& c(1+|\nabla Z_n|^2+|\frac{\partial Z_n}{\partial z}|^2+|\frac{\partial Z_n}{\partial z}|^4\|Z_n\|^4+|\frac{\partial Z_n}{\partial z}|^4|\tau|^4_{L^2(D)})| h_n|^2\\ \notag
 &&+c|\nabla Z_n|^2(1+|Z_n|^2+|\tau|^2_{L^2(D)}).
\end{eqnarray}
Let
\[
\alpha(s)=|h_n(0)|^2+c\int^s_0|\nabla Z_n|^2(1+|Z_n|^2+|\tau|^2_{L^2(D)})dr,
\]
\[
\beta(s)=c(1+|\nabla Z_n|^2+|\frac{\partial Z_n}{\partial z}|^2+|\frac{\partial Z_n}{\partial z}|^4\|Z_n\|^4+|\frac{\partial Z_n}{\partial z}|^4|\tau|^4_{L^2(D)}),
\]
applying the Gronwall inequality for (\ref{eq-23}), we conclude that
\[
|h_n(t)|^2\leq \alpha(t)+\int^t_{t_0}\alpha(s)\beta(s)\exp(\int^t_s\beta(r)dr)ds, \quad t_0\leq t\leq \Upsilon^*.
\]
By Lemma \ref{Guo}, the right hand side is uniformly bounded in $n$ as $t\rightarrow \Upsilon^*$, we conclude that $h_n$ must exit globally, i.e., $\Upsilon^*=+\infty$. Therefore, for any $\Upsilon>t_0$, we have
\begin{eqnarray}\label{eq-27}
|h_n(t)|^2\leq K_1(\Upsilon, |h(0)|, |\tau|_{L^2(D)}), \quad t\in [t_0, \Upsilon].
\end{eqnarray}
By integrating (\ref{eq-23}) with respect to $t$ over $[t_0,\Upsilon]$ and by (\ref{eq-27}), we have
\begin{eqnarray}\label{eq-28}
\int^{\Upsilon}_{t_0}\|h_n(s)\|^2ds\leq K_2(\Upsilon,|h(0)|,|\tau|_{L^2(D)}).
\end{eqnarray}
Thus, we have that $h_n$ exists globally in time and is uniformly bounded in $n$, in $L^{\infty}(t_0, \Upsilon; H)\cap L^2(t_0, \Upsilon; V)$.
\textbf{\emph{Step 3. A prior estimates for $\frac{\partial h_n}{\partial t}$.}}
Now, in order to obtain the continuity of $h_n$, we only need to show that $\partial_t h_n$ is uniformly bounded in $n$, in the space $L^2(t_0, \Upsilon; H')$, here $H'$ is the dual space of $H$. From (\ref{eq-24}), we have, for every $\psi\in C^{\infty}(\mathcal{O})$,
\[
\langle\frac{\partial h_n}{\partial t},\psi\rangle= \langle\alpha Z_n,\psi\rangle-\langle L_2 h_n,\psi\rangle-\Big\langle [(v_n\cdot \nabla)(h_n+Z_n)+\Phi(v_n)\frac{\partial (h_n+Z_n)}{\partial z}],\psi\Big\rangle,
\]
where $\langle\cdot,\cdot\rangle$ is the dual norm of $H'$ and $H$.

By Cauchy-Schwarz inequality and $\|\psi\|\geq \lambda_1 |\psi|$, we have
\begin{eqnarray}\label{eq-58}
|\langle\alpha Z_n,\psi\rangle|\leq c|Z_n||\psi|\leq c|Z_n|\|\psi\|,
\end{eqnarray}
 and by integration by parts, we have
 \begin{eqnarray}\label{eq-59}
 |\langle L_2 h_n,\psi\rangle|\leq \|h_n\|\|\psi\|.
 \end{eqnarray}
It is easy to know
\begin{eqnarray*}
&&|\Big\langle [(v_n\cdot \nabla)(h_n+Z_n)+\Phi(v_n)\frac{\partial (h_n+Z_n)}{\partial z}],\psi\Big\rangle|\\
&=&|\langle[(v_n\cdot \nabla)(h_n+Z_n)+\Phi(v_n)\frac{\partial (h_n+Z_n)}{\partial z}],\psi_n\rangle|\\
&=&|\langle[(v_n\cdot \nabla)\psi_n+\Phi(v_n)\frac{\partial \psi_n}{\partial z}],(h_n+Z_n)\rangle|,
\end{eqnarray*}
where $\psi_n=P_n \psi$.
By H\"{o}lder inequality, we have
\begin{eqnarray*}
|\int_{\mathcal{O}}(v_n\cdot \nabla)\psi_n (h_n+Z_n)dxdydz|
&\leq& c\|\psi_n\||h_n+Z_n|_3|v_n|_6\\
&\leq& c\|\psi_n\|(|h_n|^{\frac{1}{2}}\|h_n\|^{\frac{1}{2}}+|Z_n|^{\frac{1}{2}}\|Z_n\|^{\frac{1}{2}})
(|h_n|+|Z_n|+|\tau|_{L^2(D)}),
\end{eqnarray*}
and by Lemma \ref{Guo}, (\ref{eq-80}) and (\ref{eq-81}), we obtain
\begin{eqnarray*}
|\int_{\mathcal{O}}\Phi(v_n)\frac{\partial \psi_n}{\partial z}(h_n+Z_n)dxdydz|
&\leq& c|\frac{\partial \psi_n}{\partial z}||\nabla v_n|\|v_n\|^{\frac{1}{2}}_{H^2(\mathcal{O})}|h_n+Z_n|^{\frac{1}{2}}\|h_n+Z_n\|^{\frac{1}{2}}\\
&\leq& c|\frac{\partial \psi_n}{\partial z}|(|h_n|+|Z_n|+|\tau|_{L^2(D)})(\|h_n\|+\|Z_n\|+|\tau|_{L^2(D)}),
\end{eqnarray*}
then, we have
\begin{eqnarray*}
&&|\langle [(v_n\cdot \nabla)(h_n+Z_n)+\Phi(v_n)\frac{\partial (h_n+Z_n)}{\partial Z_n}],\psi\rangle|\\
&\leq &c(|h_n|+|Z_n|+|\tau|_{L^2(D)})(\|h_n\|+\|Z_n\|+|\tau|_{L^2(D)}
+|h_n|^{\frac{1}{2}}\|h_n\|^{\frac{1}{2}}+|Z_n|^{\frac{1}{2}}\|Z_n\|^{\frac{1}{2}})\|\psi_n\|,
\end{eqnarray*}
since $\psi\in H^1(\mathcal{O})$, refer to \cite{C-T-2}, we get
\[
\|\psi_n\|\leq c\|\psi\|,
\]
thus,
\begin{eqnarray}\label{eq-60}
&&|\langle [(v_n\cdot \nabla)(h_n+Z_n)+\Phi(v_n)\frac{\partial (h_n+Z_n)}{\partial z}],\psi\rangle|\\ \notag
&\leq& c(|h_n|+|Z_n|+|\tau|_{L^2(D)})(\|h_n\|+\|Z_n\|+|\tau|_{L^2(D)}
+|h_n|^{\frac{1}{2}}\|h_n\|^{\frac{1}{2}}+|Z_n|^{\frac{1}{2}}\|Z_n\|^{\frac{1}{2}})\|\psi\|.
\end{eqnarray}
Thanks to (\ref{eq-58}), (\ref{eq-59}) and (\ref{eq-60}), we have
\begin{eqnarray*}
&&|\langle\partial_t h_n,\psi\rangle|\\
&\leq& c|Z_n|\|\psi\|+c\|h_n\|\|\psi\|+c(|h_n|+|Z_n|+|\tau|_{L^2(D)})(1+\|h_n\|+\|Z_n\|+|\tau|_{L^2(D)})\|\psi\|,
\end{eqnarray*}
due to (\ref{eq-27}) and (\ref{eq-28}), we get
\begin{eqnarray}
\int^{\Upsilon}_{t_0}\|\partial_t h_n(t)\|^2_{H'}dt\leq K_3(\Upsilon,|h(0)|,|\tau|_{L^2(D)}),
\end{eqnarray}
where
\begin{eqnarray*}
K_3(\Upsilon,|h(0)|,|\tau|_{L^2(D)})&=&CK_2+c(\Upsilon-t_0)\Big(\sup_{t\in [t_0,\Upsilon]}|Z(t)|\Big)\\
&&+c\Big[K_1+\sup_{t\in [t_0,\Upsilon]}|Z(t)|+|\tau|_{L^2(D)}\Big]\Big[K_2+(1+\sup_{t\in [t_0,\Upsilon]}\|Z\|+|\tau|_{L^2(D)})(\Upsilon-t_0)\Big].
\end{eqnarray*}
Thus, we conclude that $\partial_t h_n$ is uniformly bounded in $n$ in $L^2(t_0,\Upsilon; H')$.

\textbf{\emph{Step 4. Passage to the limit $n\rightarrow \infty$.}}
 by Alaoglu compactness theorem, we can extract a subsequence $(v_{n'},T_{n'} )$ such that
\begin{eqnarray*}
T_{n'}&\rightharpoonup& T \ weakly\ in \ L^2(t_0,\Upsilon; V) \ and \ weak-star \ in\  L^{\infty}(t_0,\Upsilon; H),\\
v_{n'}&\rightharpoonup& v\ weakly\ in\ L^2(t_0,\Upsilon; H^2(\mathcal{O})) \ and\ weak-star \ in\  L^{\infty}(t_0,\Upsilon; H^1(\mathcal{O})),
\end{eqnarray*}
where
\begin{eqnarray*}
T&\in& L^{\infty}(t_0,\Upsilon; H)\bigcap L^2(t_0,\Upsilon; V),\\
v&\in &L^{\infty}(t_0,\Upsilon; H^1(\mathcal{O}))\bigcap L^{2}(t_0,\Upsilon; H^2(\mathcal{O})).
\end{eqnarray*}
From \emph{Step 3} and using the standard argument used in proving global well-posedness of 2D stochastic Navier-Stokes equations,
we can prove that $(v,T)$ is a required weak solution of
(\ref{eq-35})-(\ref{equ-12}), that is, $(p_s, v, T)$ is a weak solution of the system (\ref{eq-8})-(\ref{eq-14}).

\textbf{(Uniqueness of weak solution)}  Suppose $(h^{(i)},v^{(i)},p_s^{(i)})$, $i=1,2$ are two weak solutions of (\ref{eq-35})-(\ref{equ-12}), with the initial data $h^{(i)}_0\in H$, let
\[
\iota=h^{(1)}-h^{(2)},\ \kappa=v^{(1)}-v^{(2)},\ p_s={p_s}^{(1)}-{p_s}^{(2)},
\]
then we have
\begin{eqnarray}\label{equ-16}
\frac{d\iota}{dt}+L_2\iota+(\kappa\cdot\nabla)(h^{(1)}+Z_n)+(v^{(2)}\cdot \nabla)\iota+\Phi(\kappa)\frac{\partial(h^{(1)}+Z_n)}{\partial z}+\Phi(v^{(2)})\frac{\partial \iota}{\partial z}&=&0,\\
\label{equ-17}
\nabla\Big[p_s(x,y,t)-\int^{z}_{-1} \iota dz'\Big]+fk\times \kappa +L_1 \kappa&=&0,\\
\label{equ-18}
\frac{\partial \kappa}{\partial z}\Big|_{z=0}=0, \quad  \frac{\partial \kappa}{\partial z}\Big|_{z=-1}=0, \quad \kappa\cdot n\Big|_{\Gamma_l}=0,\quad \frac{\partial \kappa}{\partial n}\times n\Big|_{\Gamma_l}&=&0,
\end{eqnarray}
by similar method as Proposition \ref{prp-2}, we obtain
 \begin{eqnarray}\label{equa-2}
|\kappa|_{H^1(\mathcal{O})}\leq c|\iota|,\quad |\kappa|_{H^2(\mathcal{O})}\leq c\|\iota\|.
 \end{eqnarray}
 By taking inner product of (\ref{equ-16}) with $\iota$ in $L^2(\mathcal{O})$, we have
\begin{eqnarray*}
\frac{1}{2}\frac{d|\iota|^2}{dt}+\|\iota\|^2=-\Big\langle(\kappa\cdot\nabla)(h^{(1)}+Z_n),\iota\Big\rangle-\Big\langle(v^{(2)}\cdot \nabla)\iota,\iota\Big\rangle-\Big\langle\Phi(\kappa)\frac{\partial(h^{(1)}+Z_n)}{\partial z},\iota\Big\rangle-\Big\langle\Phi(v^{(2)})\frac{\partial \iota}{\partial z},\iota\Big\rangle.
\end{eqnarray*}
By integration by parts and (\ref{equ-18}), we reach
\[
\Big\langle(v^{(2)}\cdot \nabla)\iota,\iota\Big\rangle+\Big\langle\Phi(v^{(2)})\frac{\partial \iota}{\partial z},\iota\Big\rangle=0,
\]
further, applying the H\"{o}lder inequality, (\ref{equa-2}) and the Young inequality, we obtain
\begin{eqnarray*}
\left|\Big\langle(\kappa\cdot\nabla)(h^{(1)}+Z_n),\iota\Big\rangle\right|
&\leq& c\|h^{(1)}+Z_n\||\kappa|_{H^1(\mathcal{O})}\|\iota\|\\
&\leq& c|\iota|\|h^{(1)}+Z_n\|\|\iota\|\\
&\leq& \varepsilon\|\iota\|^2+c\|h^{(1)}+Z_n\|^2|\iota|^2,
\end{eqnarray*}
similarly,
\begin{eqnarray*}
\left|\Big\langle\Phi(\kappa)\frac{\partial(h^{(1)}+Z_n)}{\partial z},\iota\Big\rangle\right|
&\leq& c|\frac{\partial(h^{(1)}+Z_n)}{\partial z}||\nabla \kappa|^{\frac{1}{2}}|\kappa|^{\frac{1}{2}}_{H^2(\mathcal{O})}|\iota|^{\frac{1}{2}}\|\iota\|^{\frac{1}{2}}\\
&\leq& c| \iota||\frac{\partial(h^{(1)}+Z_n)}{\partial z}|\|\iota\|\\
&\leq& \varepsilon\|\iota\|^2+c|\frac{\partial(h^{(1)}+Z_n)}{\partial z}|^2|\iota|^2,
\end{eqnarray*}
then
\begin{eqnarray}\label{equ-14}
\frac{d|\iota|^2}{dt}+\|\iota\|^2\leq cg(t)|\iota|^2,
\end{eqnarray}
where
\[
g(t)=\|h^{(1)}+Z_n\|^2+|\frac{\partial(h^{(1)}+Z_n)}{\partial z}|^2.
\]
Notice that for the weak solution $h^{(i)}$, $i=1,2$ and $Z_n\in {D}(A)$, we have
\[
\int^t_{t_0}g(s)ds<\infty,
\]
then, by a direct integration of (\ref{equ-14}), we have
\[
|\iota(t)|^2\leq|\iota(0)|^2e^{c\int^t_{t_0}g(s)ds},
\]
thus, the uniqueness and continuity on the initial data of $\iota$ are obtained. By (\ref{equa-2}), we have the uniqueness of $v$, then (\ref{equ-10}) implies the uniqueness of $p_s$.
$\hfill\blacksquare$
\begin{thm}\label{thm-5}(Global well-posedness of strong solution)
Suppose that $\tau\in L^2(D)$. Then for every $T_0\in V$, there is a unique strong solution $(p_s,v,T)$ of the system (\ref{eq-8})-(\ref{eq-14}) in the sense of definition \ref{dfn}.
\end{thm}
\textbf{Proof of theorem \ref{thm-5}} \quad Suppose that $(p_s,v,h)$ is the weak solution of (\ref{eq-35})-(\ref{equ-12}) with the initial data $h(0)\in V$.
Taking the $L^2(\mathcal{O})$ inner product of equation (\ref{eq-35}) with $L_2 h$, we have
\[
\frac{1}{2}\frac{d\|h\|^2}{dt}+|L_2 h|^2=\int_{\mathcal{O}} \alpha Z L_2 h dxdydz-\int_{\mathcal{O}} (v\cdot \nabla)(h+Z) L_2 h dxdydz-\int_{\mathcal{O}} \Phi(v)\frac{\partial (h+Z)}{\partial z} L_2 h dxdydz.
\]
From Cauchy-Schwarz inequality,
\[
|\int_{\mathcal{O}} \alpha Z L_2 h dxdydz|\leq C|Z||L_2 h|.
\]
By (\ref{eq-34}) and  (\ref{eq-80}), we have
\begin{eqnarray*}
|\int_{\mathcal{O}}(v\cdot \nabla)h L_2hdxdydz|
&\leq& c|L_2h||\nabla h|_3|v|_6\\
&\leq& c|L_2h||\nabla h|^{\frac{1}{2}}|L_2 h|^{\frac{1}{2}}(|h|+|Z|+|\tau|_{L^2(D)})\\
&\leq& \varepsilon|L_2h|^2+c\|h\|^2(|h|^4+|Z|^4+|\tau|^4_{L^2(D)}),
\end{eqnarray*}
and similarly,
\begin{eqnarray*}
|\int_{\mathcal{O}}(v\cdot \nabla)Z L_2hdxdydz|
&\leq& c|L_2h||\nabla Z|_3|v|_6\\
&\leq& c|L_2h||\nabla Z|^{\frac{1}{2}}|Z|^{\frac{1}{2}}_{H^2(\mathcal{O})}(|h|+|Z|+|\tau|_{L^2(D)})\\
&\leq& \varepsilon|L_2h|^2+c|\nabla Z||Z|_{H^2(\mathcal{O})}(|h|^2+|Z|^2+|\tau|^2_{L^2(D)}).
\end{eqnarray*}
Applying (\ref{eq-37}) and (\ref{eq-80}), we get
\begin{eqnarray*}
|\int_{\mathcal{O}} \Phi(v)\frac{\partial h}{\partial z} L_2 h dxdydz|
&\leq& c|L_2h||\frac{\partial h}{\partial z}|^{\frac{1}{2}}|\nabla\frac{\partial h}{\partial z}|^{\frac{1}{2}}|\nabla v|^{\frac{1}{2}}|v|^{\frac{1}{2}}_{H^2(\mathcal{O})}\\
&\leq& c|L_2h|^{\frac{3}{2}}\|h\|^{\frac{1}{2}}(|h|^{\frac{1}{2}}\|h\|^{\frac{1}{2}}+\|Z\|+|\tau|_{L^2(D)})\\
&\leq& \varepsilon|L_2h|^2+c(|h|^2\|h\|^2+\|Z\|^4+|\tau|^4_{L^2(D)})\|h\|^2,
\end{eqnarray*}
and similarly,
\begin{eqnarray*}
|\int_{\mathcal{O}} \Phi(v)\frac{\partial Z}{\partial z} L_2 h dxdydz|
&\leq& c|L_2h||\frac{\partial Z}{\partial z}|^{\frac{1}{2}}|\nabla\frac{\partial Z}{\partial z}|^{\frac{1}{2}}|\nabla v|^{\frac{1}{2}}|v|^{\frac{1}{2}}_{H^2(\mathcal{O})}\\
&\leq& c|L_2h||\frac{\partial Z}{\partial z}|^{\frac{1}{2}}|\nabla\frac{\partial Z}{\partial z}|^{\frac{1}{2}}(|h|^{\frac{1}{2}}\|h\|^{\frac{1}{2}}+\|Z\|+|\tau|_{L^2(D)})\\
&\leq& \varepsilon|L_2h|^2+c|\frac{\partial Z}{\partial z}|^2|\nabla\frac{\partial Z}{\partial z}|^2\|h\|^2+c|h|^2+c|\frac{\partial Z}{\partial z}||\nabla\frac{\partial Z}{\partial z}|(\|Z\|^2+|\tau|^2_{L^2(D)}).
\end{eqnarray*}
Thus, from all the estimates, we obtain
\begin{eqnarray*}
\frac{d\|h\|^2}{dt}+|L_2 h|^2&\leq& c(1+|h|^2\|h\|^2+\|Z\|^4+|\tau|^4_{L^2(D)}+|\frac{\partial Z}{\partial z}|^2|\nabla\frac{\partial Z}{\partial z} |^2)\|h\|^2\\
&& +c(1+|\nabla Z||AZ|)|h|^2+c|\nabla Z||AZ|(1+\|Z\|^2+|\tau|^2_{L^2(D)}),
\end{eqnarray*}
from Gronwall inequality and Theorem \ref{thm-4}, for $t\in [t_0, \Upsilon]$,
\begin{eqnarray}\label{eq-38}
\|h(t)\|^2+\int^t_{t_0}| L_2h(s)|^2ds\leq K_4(\Upsilon, h(0),|\tau|_{L^2(D)}),
\end{eqnarray}
where
\begin{eqnarray*}
K_3(\Upsilon, h(0),|\tau|_{L^2(D)})&=&\Big[K_1+c(\Upsilon-t_0)\sup_{t\in [t_0, \Upsilon]}(1+|\nabla Z||AZ|)K_1+c\sup_{t\in [t_0, \Upsilon]}(|\nabla Z||AZ|)(1+\|Z\|^2+|\tau|^2_{L^2(D)})\Big]\\
&&\exp\left(c(\Upsilon-t_0)(1+\sup_{t\in [t_0, \Upsilon]}(\|Z\|^4+|\frac{\partial Z}{\partial z}|^2|\nabla\frac{\partial Z}{\partial z} |^2)+|\tau|^4_{L^2(D)})+K_1K_2\right),
\end{eqnarray*}
and $K_1$, $K_2$ are defined as (\ref{eq-27}) and (\ref{eq-28}) respectively. In addition, using similar argument as in Theorem \ref{thm-4}, we can show that
\[
\partial_z h\in L^2(t_0,\Upsilon;H).
\]
Therefore, $h\in C(t_0,\Upsilon; V)\bigcap L^2(t_0,\Upsilon; H^2(\mathcal{O}))$. Moreover, by (\ref{eq-80}), we have
\[
v\in C(t_0,\Upsilon;V)\bigcap L^2(t_0,\Upsilon;H^2(\mathcal{O})).
\]
That is, $(v,T)$ is a strong solution. Since the strong solution must be a weak solution, by Theorem \ref{thm-4}, there is only one weak solution, we conclude that the strong solution is unique.
$\hfill\blacksquare$
\section{Global Attractor}
\subsection{Preliminaries}

Let $\Big((\Omega,\mathcal{F},\mathbb{P}),(\theta_t)_{t\in \mathbb{R}}\Big)$ be a metric dynamical system over a complete probability space $(\Omega,\mathcal{F},\mathbb{P})$, i.e. $(t,\omega)\rightarrow \theta_t(\omega)$ is $\mathcal{B}(\mathbb{R})\otimes \mathcal{F}/\mathcal{F}$-measurable, $\theta_0=id$ and $\theta_{t+s}=\theta_t\circ \theta_s$ for all $t,s \in \mathbb{R}$. We recal some concepts of random dynamical systems, the detail can be referring to \cite{RDS} and \cite{CDF} ets.

\begin{dfn}\label{dfn-1}
Let (X,d) be a complete separable metric space. A family of maps $S(t,s;\omega): X\rightarrow X$, $s\leq t$ is said to be a stochastic flow, if for all $\omega\in \Omega$\begin{description}
                        \item[(i)] $S(s,s;\omega)=id$, \ for \ all \ $s\in \mathbb{R}$,
                        \item[(ii)] $S(t,s;\omega)x=S(t,r;\omega)S(r,s;\omega)x$, \ for \ all \ $t\geq r\geq s$, \ $x\in X$.
                      \end{description}
\end{dfn}
\begin{dfn}
Let (X,d) be as Definition \ref{dfn-1}. A random dynamical system (RDS) over $\theta_t$ is a measurable map
\[
\varphi : \mathbb{R}^{+} \times X \times \Omega \rightarrow X, (t,x, \omega)\rightarrow \varphi(t,\omega)x,
\]
such that $\varphi(0,\omega)=id $ and\ $\varphi$ satisfies the cocycle property, i.e.
\[
\varphi(t+s,\omega)=\varphi(t,\theta_s \omega)\circ \varphi(s,\omega),
\]
for all $t, s \in \mathbb{R}^+$ and all $\omega\in \Omega$. $\varphi$ is said to be a continuous RDS if $\mathbb{P}-$ a.s. $x\rightarrow \varphi(t,\omega)x$ is continuous for all $t\in \mathbb{R}^+$.
\end{dfn}
With the notion of an RDS at our disposal, we can recall the stochastic generalization of notions of absorption, attraction and $\Omega-$ limit sets.
\begin{dfn}
Let $(X,d)$ be as Definition \ref{dfn-1},
\begin{description}
  \item[(i)] A set-valued map $K:\Omega\rightarrow 2^{X}$ is called measurable if for all $x\in X$ the map $\omega\rightarrow d(x,K(\omega))$ is measurable, where for nonempty sets $A,B\in 2^X$,
      \[
      d(A,B)=\sup_{x\in A} \inf_{y\in B}d(x,y) \ and \ d(x,B)=d(\{x\},B).
      \]
      A measurable set-valued map is also called a random set.
  \item[(ii)] Let $A,B$ be random sets. $A$ is said to absorb $B$ if $\mathbb{P}-$ a.s. there exists an absorption time $t_B(\omega)\geq 0$ such that for all $t\geq t_B(\omega)$
      \[
      \varphi(t,\theta_{-t}\omega)B(\theta_{-t}\omega)\subseteq A(\omega).
      \]
      A is said to attract B if
      \[
      d(\varphi(t,\theta_{-t}\omega)B(\theta_{-t}\omega),A(\omega))\rightarrow 0, \ as \ t\rightarrow \infty,\ \mathbb{P}-a.s.
      \]
  \item[(iii)] For a random set $A$ we define the $\Omega-$ limit set to be
  \[
  \Omega_A(\omega)\triangleq \bigcap_{T\geq 0}\overline{\bigcup_{t\geq T}\varphi(t,\theta_{-t}\omega)A(\theta_{-t}\omega)}.
  \]
\end{description}

\end{dfn}

\begin{dfn}
A random attractor for an RDS $\varphi$ is a compact random set $\mathcal{A}$ satisfying $\mathbb{P}-$a.s.
\begin{description}
  \item[(i)] $\mathcal{A}$ is invariant, i.e. $\varphi(t,\omega)\mathcal{A}(\omega)=\mathcal{A}(\theta_t \omega)$,
  \item[(ii)] $\mathcal{A}$ attracts all nonrandom bounded sets $B\subset X$.
\end{description}
\end{dfn}
Referring to \cite{CDF}, there has  the sufficient condition for the existence of global attractor.
\begin{thm}\label{thm-6}\cite{CDF}
Let $\varphi$ be a continuous RDS and assume that there exists a compact random set $K$ absorbing all nonrandom bounded sets $B\subset X$. We set
\[
\mathcal{A}(\omega)=\overline{\bigcup_{B\subset X, \ B\ bounded}\Omega_B(\omega)},
\]
where the union is taken over all the bounded subsets of $X$.
Then for $\mathbb{P}$-a.s.
\begin{description}
  \item[(1)] $\mathcal{A}(\omega)$ is an nonempty compact invariant subset of X. In particular, if X is connected, then it is a connected subset of $K(\omega)$.
  \item[(2)] The family $\mathcal{A}(\omega)$ , $\omega \in \Omega $ is measurable.

\item[(3)] It is the minimal closed set attracting all nonrandom bounded set $B\subset X$.
\end{description}
\end{thm}
\subsection{Existence of Global Attractor}
Define the Wiener shift
\begin{eqnarray}\label{eq-99}
\theta_t\omega(s)=\omega(t+s)-\omega(t) \quad t,s\in \mathbb{R},
\end{eqnarray}
then, $\theta_s: \Omega \rightarrow \Omega $ is measure preserving and ergodic operators with respect to $\mathbb{P}$.

\textbf{Proof of Theorem \ref{thm-2}} \quad
 By Theorem \ref{thm-1}, we know that if $s\in \mathbb{R}$ and the initial data $h_s\in V$ a.s., there exists a unique strong solution $h(t,\omega)$ defined on $[s,\infty)$ to (\ref{eq-35})-(\ref{equ-12}) such that
\begin{eqnarray}\label{eq-36}
h(s, \omega)=h_s(\omega) \quad \mathbb{P}-a.s.,
\end{eqnarray}
and $h(t,\omega)$ depends continuously on the initial data $h_s$. Define the stochastic dynamical system $(S(t,s; \omega))_{t\geq s,\omega\in \Omega}$ by
\begin{eqnarray*}
S(t,s;\omega)T_s&=&h(t,\omega)+Z(t,\omega),
\end{eqnarray*}
where $h(t,\omega)$ is the solution of (\ref{eq-35})-(\ref{equ-12}), and
\begin{eqnarray*}
Z(t,\omega)=\int^t_{-\infty}e^{(t-u)(-A-\alpha)}GdW(u).
\end{eqnarray*}
By the uniqueness of $h(t,\omega)$, we have for all $\omega\in \Omega$, $r,s,t\in \mathbb{R}$, $s\leq r\leq t$,
\begin{eqnarray}\notag
S(s,s;\omega)&=&i_d,\\
S(t,s;\omega)&=&S(t,r;\omega)S(r,s;\omega),\label{eq-97}\\
S(t,s;\omega)&=&S(t-s,0;\theta_s\omega),\label{eq-96}
\end{eqnarray}
that is, $(S(t,s; \omega))_{t\geq s,\omega\in \Omega}$ is a stochastic flow. Furthermore, by (\ref{eq-97}) and (\ref{eq-96}), for any $s,t \in \mathbb{R}^{+}$, $T_0\in V$, we have $\mathbb{P}-a.s.$
\begin{eqnarray}\label{eq-98}
S(t+s,0;\omega)T_0&=&S(t,0;\theta_s\omega)S(s,0;\omega)T_0.
\end{eqnarray}
Define $\varphi:\mathbb{R}^+\times V\times \Omega\rightarrow V$,
\[
\varphi(t,\omega)T_0=S(t,0;\omega)T_0,
\]
then, (\ref{eq-98}) implies the cocycle property $\varphi(t+s,\omega)=\varphi(t,\theta_s\omega)\circ \varphi(s,\omega)$ holds for $\varphi$.
The joint measurability of $\varphi$ follows from the construction of the solutions $h(t)$. Hence $\varphi$ is a continuous RDS. Based on Theorem \ref{thm-6}, we only need to find a compact attracting set $K(\omega)$ at time 0.

Let $B$ be a bounded set in $V$ and $T_s\in B$ for any $s\in \mathbb{R}$. Let $h(t,\omega)$ be the solution of (\ref{eq-35})-(\ref{equ-12}) with the initial data $h_s=T_s-Z(s,\omega)$.
Fix $\omega\in \Omega$, we multiply (\ref{eq-35}) by $h$ in $H$, we obtain
\[
\frac{1}{2}\frac{d|h|^2}{dt}+\|h\|^2+\left\langle\Big((v\cdot\nabla)(h+Z)+\Phi(v)\frac{\partial(h+Z)}{\partial z}\Big),h\right\rangle=\alpha\langle Z,h\rangle,
\]
by integration by parts,
\[
\left\langle\Big((v\cdot\nabla)h+\Phi(v)\frac{\partial h}{\partial z}\Big),h\right\rangle=0,
\]
then using (\ref{eq-34})and the Young inequality, we obtain
\begin{eqnarray*}
\left|\Big\langle(v\cdot\nabla)Z,h\Big\rangle\right|&\leq& c|h||AZ||v|_{H^1(\mathcal{O})}
\leq c|h|^2|AZ|^2+c|v|_{H^1(\mathcal{O})}^2,
\end{eqnarray*}
furthermore, by (\ref{eq-37}), we get
\begin{eqnarray*}
\left|\Big\langle\Phi(v)\frac{\partial Z}{\partial z},h\Big\rangle\right|&\leq& c|h||AZ||v|_{H^1(\mathcal{O})}^{\frac{1}{2}}|v|^{\frac{1}{2}}_{H^2(\mathcal{O})}
\leq c|h||AZ||v|_{H^1(\mathcal{O})}^{\frac{1}{2}}\|h+Z\|^{\frac{1}{2}}\\
&\leq& \varepsilon\|h\|^2+c|h|^2|AZ|^2+c|v|_{H^1(\mathcal{O})}^2+c\|Z\||v|_{H^1(\mathcal{O})},
\end{eqnarray*}
as the strong solution $v\in C(s,\infty; \breve{V})$, it follows
\begin{eqnarray}\label{equation1}
\frac{d|h|^2}{dt}+\|h\|^2\leq c|AZ|^2|h|^2+\frac{\lambda_1}{4}|h|^2+g,
\end{eqnarray}
with
\[
g=\frac{4\alpha^2}{\lambda_1}|Z|^2+c\|Z\|^2,
\]
where $\lambda_1$ is the first eigenvalue of $A$, which satisfies
\[
\|h\|^2\geq \lambda_1|h|^2\quad h\in V.
\]
From (\ref{equation1}), we deduce that
\begin{equation}\label{equ-5}
\frac{d|h|^2}{dt}+\frac{1}{2}\|h\|^2+\Big(\frac{\lambda_1}{4}-2c|AZ|^2\Big)|h|^2\leq g,
\end{equation}
and by Gronwall lemma, for $s<-1$, $t\in [-1,0]$,
\begin{eqnarray}
|h(t)|^2
&\leq& |h(s)|^2\exp\left(-\int^{t}_{s}\Big(\frac{\lambda_1}{4}-2c|AZ(\sigma)|^2\Big)d\sigma\right)\\ \nonumber
&&\quad+\int^{t}_{s}g(\sigma)\exp\left(-\int^{t}_{\sigma}(\frac{\lambda_1}{4}-2c|AZ(\tau)|^2)d\tau\right)d\sigma \\ \nonumber
&\leq&  c|h(s)|^2\exp\left(s\Big(\frac{\lambda_1}{4}+\frac{2c}{s}\int^{0}_{s}|AZ(\sigma)|^2d\sigma\Big)\right)\\ \label{eq-29}
&&\quad+c\int^{0}_{s}g(\sigma)\exp\left(-\int^{0}_{\sigma}\Big(\frac{\lambda_1}{4}-2c|AZ(\tau)|^2\Big)d\tau\right)d\sigma. \notag
\end{eqnarray}
From the ergodicity of the process $Z$ with values in $\mathcal{D}(A)$, we have
\begin{eqnarray*}
-\frac{1}{s}\int^s_0|AZ(\sigma)|^2d\sigma\rightarrow \mathbb{E}|AZ(0)|^2,
\end{eqnarray*}
when $s\rightarrow -\infty$. Thus, there exists $s_0(\omega)<0$ such that for any $s\leq s_0(\omega)$,
\[
-\frac{1}{s}\int^s_0|AZ(\sigma)|^2d\sigma\leq 2\mathbb{E}|AZ(0)|^2,
\]
and
\begin{equation}\label{equ-15}
\exp\Big(s\big(\frac{\lambda_1}{4}+\frac{2c}{s}\int^{0}_{s}|AZ(\sigma)|^2d\sigma\big)\Big)\leq\exp\Big(s\big(\frac{\lambda_1}{4}+4c\mathbb{E}|AZ(0)|^2\big)\Big).
\end{equation}
By Proposition \ref{prp-5}, taking $\alpha$ large enough, we obtain
\begin{eqnarray}\label{eq-39}
\mathbb{E}|AZ(0)|^2\leq\frac{\lambda_1}{16c},
\end{eqnarray}
which implies that the first term of (\ref{eq-29}) decays to 0 when $s\rightarrow -\infty$.
Moreover, by laws of iterated logarithm of Brownian motion and from
\[
Z_j(t)=Z_j(0)-\alpha\int^{0}_{t}Z_j(s)ds +w_j(t),
\]
we have $\frac{|Z_j(t)|}{t}$ is bounded at $-\infty$ and that $g(t)$ grows at most polynomially.

Since $g(t)$ is multiplied by a function which delays exponentially, so by (\ref{eq-29}) and (\ref{equ-15}), the integral converges. Thus, for $s< s_0(\omega)$, $t\in [-1,0]$,
\begin{eqnarray*}
|h(t)|^2&\leq & c|h(s)|^2\exp(\frac{\lambda_1 s}{2})
+c\int^{0}_{-\infty}g(\sigma)\exp\left(\sigma\Big(\frac{\lambda_1}{4}+\frac{2c}{\sigma}\int^{0}_{\sigma}|AZ(\tau)|^2d\tau\Big)\right)d\sigma\\
&\leq & 2c|T_s|^2\exp(\frac{\lambda_1 s}{2})+2c|Z(s)|^2\exp(\frac{\lambda_1 s}{2})\\
&&\quad+c\int^{0}_{-\infty}g(\sigma)\exp\left(\sigma\Big(\frac{\lambda_1}{4}+\frac{2c}{\sigma}\int^{0}_{\sigma}|AZ(\tau)|^2d\tau \Big)\right)d\sigma.
\end{eqnarray*}
It's now clear that there exists $s_1(\omega, B)$, depending only on $B$ and $\omega$, such that for $s<s_1(\omega, B)$, $t\in[-1,0]$
\begin{eqnarray}\label{eq-40}
|h(t)|^2\leq r_0(\omega),
\end{eqnarray}
where
\begin{eqnarray*}
r_0(\omega)&=&c\int^{0}_{-\infty}g(\sigma)\exp\left(\sigma\Big(\frac{\lambda_1}{4}+\frac{2c}{\sigma}\int^{0}_{\sigma}|AZ(\tau)|^2d\tau\Big)\right)d\sigma\\
&&\quad+2c\sup_{s\in (-\infty,-1]}\Big(|Z(s)|^2\exp(\frac{\lambda_1 s}{2})\Big)+1.
\end{eqnarray*}
Moreover, we can integrate (\ref{equ-5}) about $t$ on [-1,0] to deduce
\begin{eqnarray}\label{eq-41}
\int^{0}_{-1}\|h(t)\|^2dt\leq r_1(\omega),
\end{eqnarray}
where
\[
r_1(\omega)=4c\Big(\int^0_{-1}|AZ(\sigma)|^2d\sigma\Big)r_0(\omega)+2\int^0_{-1}g(\sigma)d\sigma.
\]
To derive an estimate in $V$, we take the scalar product of (\ref{eq-35}) by $Ah$ in $H$ and obtain
\begin{equation}\notag
\frac{1}{2}\frac{d\|h\|^2}{dt}+|Ah|^2=\alpha\langle Z,Ah\rangle-\Big\langle(v\cdot\nabla)(h+Z)+\Phi(v)\frac{\partial(h+Z)}{\partial z},Ah\Big\rangle.
\end{equation}
Applying H\"{o}lder inequality and the Young inequality, we get
\[
\Big|\alpha\langle Z,Ah\rangle\Big|\leq 2\alpha^2|Z|^2+\frac{1}{8}\|h\|^2,
\]
and furthermore, by (\ref{eq-34}), we have
\begin{eqnarray*}
\left|\Big\langle(v\cdot \nabla)(h+Z),Ah\Big\rangle\right|
&\leq& c|Ah||\nabla(h+Z)|_3|v|_{6}\\
&\leq & c|Ah|\|h+Z\|^{\frac{1}{2}}|A(h+Z)|^{\frac{1}{2}}(|h+Z|+|\tau|_{L^2(D)}),
\end{eqnarray*}
by using (\ref{eq-37}), we obtain
\begin{eqnarray*}
\left|\Big\langle\Phi(v)\frac{\partial (h+Z)}{\partial z},Ah\Big\rangle\right|&\leq& c|Ah||v|_{H^1(\mathcal{O})}^{\frac{1}{2}}|v|^{\frac{1}{2}}_{H^2(\mathcal{O})}\|h+Z\|^{\frac{1}{2}}|A(h+Z)|^{\frac{1}{2}}\\
&\leq& c|Ah|(|h+Z|+|\tau|_{L^2(D)})^{\frac{1}{2}}(\|h+Z\|+|\tau|_{L^2(D)})^{\frac{1}{2}}\|h+Z\|^{\frac{1}{2}}|A(h+Z)|^{\frac{1}{2}}.
\end{eqnarray*}
As $\tau \in L^2(D)$, we omit it for written simplicity. Then
\begin{eqnarray*}
\left|\Big((v\cdot\nabla)(h+Z)+\Phi(v)\frac{\partial(h+Z)}{\partial z},Ah\Big)\right|&\leq&c|h+Z|^{\frac{1}{2}}|A(h+Z)|^{\frac{1}{2}}\|h+Z\||Ah|\\
&\leq& 2c^2|h+Z||A(h+Z)|\|h+Z\|^2+\frac{1}{8}|Ah|^2\\
&\leq& 16c^4|h+Z|^2\|h\|^4+\frac{1}{16}|Ah|^2+16c^4|h+Z|^2\|Z\|^4\\
&&\quad+\frac{1}{16}|Ah|^2+2c^2|h+Z||AZ|\|h+Z\|^2+\frac{1}{8}|Ah|^2.
\end{eqnarray*}
Collecting all the inequalities above, we deduce
\begin{eqnarray}\label{eq-42}
\frac{d\|h\|^2}{dt}+|Ah|^2\leq G(t)+H(t)\|h\|^2,
\end{eqnarray}
where
\begin{eqnarray*}
G(t)&=&4\alpha^2|Z|^2+4c^2|h+Z||AZ|\|h+Z\|^2+32c^4|h+Z|^2\|Z\|^4,\\
H(t)&=&32c^4|h+Z|^2\|h\|^2.
\end{eqnarray*}
Let
\[
K(t)=\|h(t)\|^2e^{\int^0_tH(\sigma)d\sigma}+\int^0_{t}G(\sigma)e^{\int^0_{\sigma}H(\tau)d\tau}d\sigma,
\]
by (\ref{eq-42}), we have
\[
\frac{dK(t)}{dt}\leq 0,
\]
then, for any $t\in [-1,0]$, $K(0)\leq K(t)$, that is,
\begin{eqnarray*}
\|h(0)\|^2&\leq& \|h(t)\|^2e^{\int^0_tH(\sigma)d\sigma}+\int^0_{t}G(\sigma)e^{\int^0_{\sigma}H(\tau)d\tau}d\sigma\\
&\leq&\Big(\|h(t)\|^2+\int^0_{-1}G(\sigma)d\sigma\Big)e^{\int^0_{-1}H(\sigma)d\sigma},
\end{eqnarray*}
and integrating the above inequality with respect to $t$ on [-1,0],
\begin{eqnarray*}
\|h(0)\|^2&\leq&\left(\int^0_{-1}\|h(t)\|^2dt+\int^0_{-1}G(\sigma)d\sigma\right)e^{\int^0_{-1}H(\sigma)d\sigma}.
\end{eqnarray*}
Recalling (\ref{eq-40}) and (\ref{eq-41}), then there exists $r_2(\omega)$ such that, when $s<s_1(\omega, B)$,
\begin{eqnarray}\label{equ-6}
\|h(0)\|^2\leq r_2(\omega).
\end{eqnarray}
Notice that, since $\varphi(-s,\theta_{s}\omega)=S(-s,0;\theta_{s}\omega)=S(0,s;\omega)$, thus in the language of the stochastic flow, we have
\begin{eqnarray*}
T(0,\omega)&=&S(0,s;\omega)T_s\\
&=&h(0,\omega)+Z(0,\omega).
\end{eqnarray*}
Let $K(\omega)$ be the ball in $V$ of radius $r^{\frac{1}{2}}_2(\omega)+\|Z(0,\omega)\|$, we have proved for any $B$ bounded in $V$, there exists $s_1(\omega,B)$ such that, for $s<s_1(\omega,B)$,
\[
\varphi(-s,\theta_{s}\omega)B\subset K(\omega) \quad P-a.e.
\]
This clearly implies that $K(\omega)$ is attracting at time 0.

Due to Theorem \ref{thm-6}, only the compactness is left for checking. Following the same method via a special technique using an Aubin-Lions compactness lemma combined with the use of the Riesz lemma and a continuity argument in \cite{N}, we easily obtain that $\varphi(-s,\theta_{s}\omega)$ is a compact operator from $V$ to $V$. We do not give the detail here. Thus, the proof of Theorem \ref{thm-2} is complete.
$\hfill\blacksquare$
\subsection{Hausdorff Dimension of Global Attractor}
\subsubsection{ Preliminaries}
Firstly, we recall some definitions on Hausdorff measure.

\begin{dfn}(\cite{Temam})
For a subset $Y$ of a metric space, $s\geq 0$ and $\varepsilon>0$ put
\[
\mu_h(Y,s,\varepsilon)=\inf\Big\{\sum_{i\in I}r^s_i:\ (B_i)_{i\in I}\ is\ a\ countable\ collection\ of\ balls\ of\ radii\ r_i\ covering\ Y,\ r_i\leq \varepsilon\Big\},
\]
where $\inf \emptyset=\infty$. The $s$-dimensional Hausdorff measure of $Y$ is
\[
\mu_H(Y,s)=\lim_{\varepsilon\rightarrow  0}\mu_H(Y,s,\varepsilon)=\sup_{\varepsilon >0}\mu_H(Y,s,\varepsilon).
\]
\end{dfn}
\begin{dfn} (\cite{Temam})
The Hausdorff dimension $dim _H (Y)$ of a subset $Y$ of a metric space is
\[
dim_H(Y)=\inf\{s\geq 0: \mu_H(Y,s)=0\}=\sup\{s\geq 0:\mu_H(Y,s)>0\},
\]
where $\mu_H(Y,s)$ is the $s-$ dimensional Hausdorff measure of $Y$.
\end{dfn}

We next show that volume contraction of an random dynamical system (RDS) on a random compact invariant set implies that the Hausdorff dimension of this set is finite.

Suppose $\varphi$ is a continuous RDS on a separable Hilbert space $H$ with norm $|\cdot|$, $L$ is a bounded linear operator on $H$.  $\wedge^d L$ denotes the $d-$ fold exterior product of $L$, and $\|\cdot\|_{\mathcal{L}(H^d, \wedge^d H)}$ denotes the operator norm defined by

\[
\|\wedge^d L\|_{\mathcal{L}(H^d, \wedge^d H)}=\sup_{\{\varphi_1,\cdot\cdot\cdot, \varphi_d\in H; |\varphi_i|_H\leq 1, \forall i\}}|(\wedge^d L)(\varphi_1,\cdot\cdot\cdot, \varphi_d)|_{\wedge^d H},
\]
where
\[
(\wedge^d L)(\varphi_1,\cdot\cdot\cdot, \varphi_d)=L\varphi_1\wedge\cdot\cdot\cdot L\varphi_d.
\]
Suppose that $\omega\mapsto X(\omega)$ is a compact strictly invariant set for $\varphi$.
\begin{dfn}\label{dfn-2}
(\cite{CFH})  $\varphi$ is called  weakly differentiable on $X,$ if for $\mathbb{P}-$ almost all $\omega$, every $u\in X(\omega)$ and $t>0,$ there exists a linear map $D_u \varphi(t,\omega):H\rightarrow H$ such that
\[
g_{\delta}(t,\omega)=\sup\Big\{\frac{|\varphi(t,\omega)u-\varphi(t,\omega)v-D_u \varphi(t,\omega)(u-v)|}{|u-v|}: u,v\in X(\omega), |u-v|\leq \delta\Big\}<\infty
\]
for any $\delta>0$, and converges to zero $\mathbb{P}-a.s.$ as $\delta\rightarrow 0$ for any $t>0$ fixed.
\end{dfn}
If $\varphi$ is Fr$\acute{e}$cher differentiable in $u$, then, the map $D_u \varphi(t,\omega)$ is the differential.  Put
\[
\gamma_1(t,\omega)=\sup_{u\in X(\omega)}\|D_u \varphi(t,\omega)\|, \quad \Gamma_d(t,\omega)=\sup_{u\in X(\omega)}\|\wedge^d D_u \varphi(t,\omega)\|_{\mathcal{L}(H^d, \wedge^d H)}.
\]
The following is the sufficient condition to obtain the finite Hausdorff dimension.
\begin{thm}\label{thm-7}(\cite{CFH})
Let $\varphi$ be an RDS on a separable Hilbert space $H$ and $\omega\mapsto X(\omega)$ be a random compact set, strictly invariant for $\varphi$. Suppose that $\varphi$ is weakly differentiable on X in the sense of Definition \ref{dfn-2} and there exist $d$ and $t_0>0$ such that $\mathbb{P}-$ a.s.
\[
\sup_{u\in X(\omega)}\|\wedge^d D_u \varphi(t_0,\omega)\|_{\mathcal{L}(H^d, \wedge^d H)}=\Gamma_d(t_0,\omega)<1.
\]
Suppose further that
\begin{description}
  \item[(i)] $\gamma_1(t,\omega)\in L^{\infty}(\Omega,\mathbb{P})$ for every $t\geq 0$.
  \item[(ii)] $g_{\delta}(t,\omega)\in L^{\infty}(\Omega,\mathbb{P})$ for every $t\geq 0$.
  \item[(iii)] $g_{\delta}(t,\omega)\rightarrow 0$  in $L^{\infty}(\Omega,\mathbb{P}) $ for every $t\geq 0$, as $\delta\rightarrow 0$.
\end{description}
Then, the Hausdorff dimension of $X$ satisfies $\dim_H(X)\leq d,\ \  P-a.s.$.
\end{thm}
In Sect 5.2, we have obtain the existence of global attractor $\mathcal{A}(\omega)$ in $V$, $\mathbb{P}-$ a.s.. In the following, we will use Theorem \ref{thm-7} to prove Theorem \ref{thm-8} holds.
\subsubsection{Finite Dimension of the Attractor for the System (\ref{eq-8})-(\ref{eq-14})}
Using the similar method as \cite{S-G}, we need firstly make the volume contraction for (\ref{eq-8})-(\ref{eq-14}).
\begin{prp}\label{prp-3}
Let $\varphi$ be the RDS associated with the system (\ref{eq-8})-(\ref{eq-14}), and let
$\omega\rightarrow \mathcal{A}(\omega)$ be the global attractor given by Theorem \ref{thm-2}. Then there exist a positive integer $d$ and a time $t>0$ such that
\[
|\Gamma_d(t,\cdot)|_{L^{\infty}(\Omega)}<1,
\]
where $\Gamma_d(t,\omega)=\sup_{u\in \mathcal{A}(\omega)}\|\wedge^d D_u \varphi(t,\omega)\|$.
\end{prp}
\textbf{Proof of Proposition \ref{prp-3}} \quad Consider the volume contraction for (\ref{eq-8})-(\ref{eq-14}),
\begin{eqnarray}\label{eq-43}
\nabla [q_s(x,y,t)-\int^z_{-1}\chi(x,y,z',t)dz']+fk\times u+L_1 u&=&0,\\
\label{eq-44}
\nabla \cdot \int^0_{-1} u(x,y,z',t)dz'&=&0,\\
\label{eq-45}
\partial_t \chi&=&F'(T)\chi,\\
\label{eq-45-1}
\frac{\partial u}{\partial z}|_{z=0}=0,\quad \frac{\partial u}{\partial z}|_{z=-1}=0, \quad u\cdot n|_{\Gamma_l}=0, \quad \frac{\partial u}{\partial n}\times n|_{\Gamma_l}&=&0,\\
\label{eq-46}
\chi(0,\omega)&=&\xi\in V,
\end{eqnarray}
where
\[
F'(T)\chi=-[L_2 \chi+(u\cdot)T+(v\cdot \nabla)\chi+\Phi(u)\frac{\partial T}{\partial z}+\Phi(v)\frac{\partial \chi}{\partial z}].
\]
Similar to the prove of Theorem \ref{thm-1}, we can easily obtain the existence of strong solution $(\chi, u, q_s)$,  which satisfies
 \begin{eqnarray*}
 \chi &\in& C(t_0, \infty; V)\bigcap L^2(t_0, \infty; D(A)),\\
 u&\in& C(t_0, \infty; H^1(\mathcal{O}))\bigcap L^2(t_0, \infty; H^2(\mathcal{O})),\\
 q_s &\in& C(t_0, \infty; H^1(D))\bigcap L^2(t_0, \infty; H^2(D)).
 \end{eqnarray*}
 For $\xi=\xi_1,\cdot\cdot\cdot, \xi_d\in V$, let $\chi=\chi_1,\cdot\cdot\cdot, \chi_d$ denote the corresponding solution of (\ref{eq-43})-(\ref{eq-46}), hence, following the procedure in \cite{Temam},  we consider, for any $d\in \mathbb{N}$,
 \begin{eqnarray}
 |\chi_1(t)\wedge \cdot\cdot\cdot \chi_d(t)|^2_{\wedge^d V}=|\xi_1(t)\wedge \cdot\cdot\cdot \xi_d(t)|^2_{\wedge^d V}\exp\Big(\int^t_0 Tr F'(\varphi(\tau)T_0)\circ Q_d(\tau)d\tau\Big),
 \end{eqnarray}
 where $\varphi(\tau,\omega)T_0=T(\tau,\omega)$ and $Q_d$ is the orthogonal projector in $V$ onto the space spanned by $\chi_1(\tau)\wedge \cdot\cdot\cdot \wedge\chi_d(\tau)$. At a given time $\tau$, let
 $\psi_j$, $j\in \mathbb{N}$  be an orthogonal basis of $V$. Thus, we have
\[
Tr F'(\varphi(\tau)T_0)\circ Q_d(\tau)=\sum^d_{j=1}\Big(F'(\varphi(\tau)T_0)\psi_j(\tau),\psi_j(\tau)\Big).
\]
Notice that
\begin{eqnarray*}
(F'(\varphi(\tau)T_0)\psi_j(\tau),\psi_j(\tau))&=& -|L_2 \psi_j|^2+\int_{\mathcal{O}}(\phi_j \cdot \nabla)T L_2 \psi_jdxdydz+\int_{\mathcal{O}}(v \cdot \nabla)\psi_j L_2 \psi_jdxdydz\\
&&+\int_{\mathcal{O}}\Phi(\phi_j) \frac{\partial T}{\partial z} L_2 \psi_jdxdydz+\int_{\mathcal{O}}\Phi(v) \frac{\partial\psi_j}{\partial z} L_2 \psi_jdxdydz,
\end{eqnarray*}
where, for $j=1,2,...,d$, $\phi_j(x,y,z,\tau)$ is the solution of the following linear system:
\begin{eqnarray*}
\nabla [(q_s)_j(x,y,t)-\int^z_{-1}\chi_j(x,y,z',t)dz']+fk\times u_j+L_1 u_j&=&0,\\
\nabla \cdot \int^0_{-1} u_j(x,y,z',t)dz'&=&0,\\
\frac{\partial {u_j}}{\partial z}|_{z=0}=0,\quad \frac{\partial {u_j}}{\partial z}|_{z=-1}=0, \quad u_j\cdot n|_{\Gamma_l}=0, \quad \frac{\partial {u_j}}{\partial n}\times n|_{\Gamma_l}&=&0.
\end{eqnarray*}
Here $(q_s)_j$ and $\phi_j$ are the unknowns, while $\psi_j$ is given and fixed.
Then, by Proposition 2.5 in \cite{C-T-2}, we have
 \begin{eqnarray}\label{eq-47}
 |\phi_j|_{H^1(\mathcal{O})}\leq c|\psi_j|,\quad |\phi_j|_{H^2(\mathcal{O})}\leq c\|\psi_j\|.
 \end{eqnarray}
Now, we will estimate $Tr F'(\varphi(\tau)T_0)\circ Q_d(\tau)$. By H\"{o}lder inequality and (\ref{eq-47}), we have
\[
|\int_{\mathcal{O}}(\phi_j \cdot \nabla)T L_2 \psi_jdxdydz|\leq c| L_2 \psi_j||\nabla T||\phi_j|_{\infty}\leq c| L_2 \psi_j||\nabla T|\|\psi_j\|,
\]
and
\[
|\int_{\mathcal{O}}(v \cdot \nabla)\psi_j L_2 \psi_jdxdydz|\leq c|L_2 \psi_j||\nabla \psi_j|_3|v|_6
\leq c|L_2 \psi_j|^{\frac{3}{2}}|\nabla \psi_j|^{\frac{1}{2}}\|v\|,
\]
by (\ref{eq-37}) and (\ref{eq-47}), we get
\begin{eqnarray*}
|\int_{\mathcal{O}}\Phi(\phi_j) \frac{\partial T}{\partial z} L_2 \psi_jdxdydz|
&\leq& c|L_2 \psi_j||\frac{\partial T}{\partial z}|^{\frac{1}{2}}|\nabla\frac{\partial T}{\partial z}|^{\frac{1}{2}}|\nabla \phi_j|^{\frac{1}{2}}|\phi_j|^{\frac{1}{2}}_{H^2(\mathcal{O})}\\
&\leq& c|L_2 \psi_j||\frac{\partial T}{\partial z}|^{\frac{1}{2}}|\nabla\frac{\partial T}{\partial z}|^{\frac{1}{2}}\| \psi_j\|,
\end{eqnarray*}
similarly, we obtain
\begin{eqnarray*}
|\int_{\mathcal{O}}\Phi(v) \frac{\partial\psi_j}{\partial z} L_2 \psi_jdxdydz|
&\leq& c|L_2 \psi_j||\frac{\partial \psi_j}{\partial z}|^{\frac{1}{2}}|\nabla\frac{\partial \psi_j}{\partial z}|^{\frac{1}{2}}|\nabla v|^{\frac{1}{2}}|v|^{\frac{1}{2}}_{H^2(\mathcal{O})}\\
&\leq& c|L_2 \psi_j|^{\frac{3}{2}}|\frac{\partial \psi_j}{\partial z}|^{\frac{1}{2}}|\nabla v|^{\frac{1}{2}}|v|^{\frac{1}{2}}_{H^2(\mathcal{O})}.
\end{eqnarray*}
Thus, by the Young inequality, we have
\begin{eqnarray*}
Tr F'(\varphi(\tau)T_0)\circ Q_d(\tau)
&\leq& \sum^d_{j=1}\Big[-\frac{1}{2}|L_2 \psi_j|^2+c\Big(|\nabla T|^2+\|v\|^4+| \frac{\partial T}{\partial z}||\nabla \frac{\partial T}{\partial z}|+c|\nabla v|^2|v|^2_{H^2(\mathcal{O})}\Big)\|\psi_j\|^2\Big]\\
&\leq& \Big(\sum^d_{j=1}-\frac{1}{2}|L_2 \psi_j|^2\Big)+cd\Big(|\nabla T|^2+\|v\|^4+|\frac{\partial T}{\partial z}||\nabla \frac{\partial T}{\partial z}|+|\nabla v|^2|v|^2_{H^2(\mathcal{O})}\Big),
\end{eqnarray*}
where $\|\psi_j\|=1$ is used in the second inequality and $c$ is a constant, which is independent of $d$.

Applying generalization of the Sobolev-Lieb-thirring inequality in \cite{Temam}, we have
\begin{eqnarray*}
\sum^d_{j=1}\int_{\mathcal{O}}|L_2  \psi_j|^2dxdydz
&\geq& \frac{1}{\kappa}(\int_{\mathcal{O}}(\sum^d_{j=1}|\nabla \psi_j|^2)^2dxdydz)^{\frac{2}{3}}-\frac{1}{|\mathcal{O}|^{\frac{2}{3}}}\int_{\mathcal{O}}\sum^d_{j=1}|\nabla \psi_j|^2dxdydz\\
&\geq& d^{\frac{4}{3}}|\mathcal{O}|^{\frac{2}{3}}\frac{1}{\kappa}-d|\mathcal{O}|^{\frac{1}{3}}.
\end{eqnarray*}
We conclude
\begin{eqnarray*}
|\chi_1(t)\wedge \cdot\cdot\cdot \chi_d(t)|^2_{\wedge^d V}&=&|\xi_1(t)\wedge \cdot\cdot\cdot \xi_d(t)|^2_{\wedge^d V}\exp\Big[-d^{\frac{4}{3}}|\mathcal{O}|^{\frac{2}{3}}\frac{1}{2\kappa}t+d|\mathcal{O}|^{\frac{1}{3}}t\\
&&+cd\int^t_0 \left(|\nabla T|^2+\|v\|^4+|\frac{\partial T}{\partial z}||\nabla \frac{\partial T}{\partial z}|+|\nabla v|^2|v|^2_{H^2(\mathcal{O})}\right)d\tau \Big],
\end{eqnarray*}
from Theorem \ref{thm-1},
\[
\int^t_{t_0}\big(|\nabla T|^2+\|v\|^4+|\frac{\partial T}{\partial z}||\nabla \frac{\partial T}{\partial z}|+|\nabla v|^2|v|^2_{H^2(\mathcal{O})}\big)d\tau\leq c+ct,
\]
for $\mathbb{P}-$ a.s.. Thus, for any $t$, we get
\begin{eqnarray}\label{eq-48}
\Gamma_d(t,\omega)\leq \exp\Big(-d^{\frac{4}{3}}|\mathcal{O}|^{\frac{2}{3}}\frac{1}{2\kappa}t+cd(1+t+|\mathcal{O}|^{\frac{1}{3}}t)\Big),
\end{eqnarray}
let $d$ sufficiently large, we have $|\Gamma_d(t,\omega)|_{L^{\infty}(\Omega)}<1$ for every $t>0$.

$\hfill\blacksquare$

Now, we are ready to prove Theorem \ref{thm-8}.

\textbf{Proof of Theorem \ref{thm-8}} \quad From (\ref{eq-48}) we get existence of some large $d$ such that
\[
|\Gamma_d(t,\omega)|_{L^{\infty}}<1.
\]
From the proof of (\ref{eq-48}), we get for the particular case $d=1$,
\[
\|\chi_1\|\leq \|\xi_1\|\exp (ct+c),
\]
thus, we have $\gamma_1(t)\in L^{\infty}(\Omega, \mathbb{P})$. Finally, we have to prove that $g_{\delta}(t)$ converges to 0 in $L^{\infty}(\Omega, \mathbb{P})$ for each $t>t_0$ as $\delta\rightarrow 0$.

For $\xi \in V$, $\xi_0 \in V$, let
\begin{eqnarray*}
\chi (t)&=&\varphi(t,\omega)\xi, \quad \chi_0 (t)=\varphi(t,\omega)\xi_0,\\
 \theta&=& D_{\xi_0}\chi, \quad \theta_0= D_{\xi_0}\chi_0,\quad \bar{\theta}=\theta-\theta_0,\\
 \delta&=&D_{\xi_0}v, \quad \delta_0=D_{\xi_0}v_0,\quad \bar{\delta}=\delta-\delta_0,\\
\eta(t)&=&\chi-\chi_0-\bar{\theta},\quad \gamma(t)=v-v_0-\bar{\delta},
\end{eqnarray*}
from (\ref{eq-43})-(\ref{eq-46}), we obtain
\begin{eqnarray*}
-\int^z_{-1} \nabla \theta dz'+fk\times \delta +L_1 \delta&=&0,\\
\int^z_{-1}\nabla \cdot \delta dz'&=&0,\\
\frac{d\theta}{d t}+(\delta \cdot \nabla) \chi+(v\cdot \nabla)\theta+\Phi(\delta)\frac{\partial \chi}{\partial z}+\Phi(v)\frac{\partial \theta}{\partial z}+L_2 \theta&=&0,\\
\frac{\partial \delta}{\partial z}|_{z=0}=0,\quad \frac{\partial \delta}{\partial z}|_{z=-1}=0, \quad \delta\cdot n|_{\Gamma_l}=0, \quad \frac{\partial \delta}{\partial n}\times n|_{\Gamma_l}&=&0,\\
\frac{\partial \theta}{\partial z}|_{z=0}=0,\quad \frac{\partial \theta}{\partial z}|_{z=-1}=0, \quad \frac{\partial \theta}{\partial n}|_{\Gamma_l}&=&0,
\end{eqnarray*}
with $\theta(0)=\xi_0\in V$, using the similar method as Theorem \ref{thm-1}, we get
\begin{eqnarray}\label{eq-51}
\delta &\in& C(t_0, \Upsilon; H^1(\mathcal{O}))\bigcap L^2(t_0, \Upsilon; H^2(\mathcal{O})),\\
\label{eq-52}
\theta &\in& C(t_0, \Upsilon; V)\bigcap L^2(t_0, \Upsilon; H^2(\mathcal{O})).
\end{eqnarray}
Further, from (\ref{eq-8})-(\ref{eq-14}), it gives
\begin{eqnarray*}
\nabla \Big(q_s-\int^z_{-1}(\chi-\chi_0)dz'\Big)+fk\times (v-v_0)+L_1(v-v_0)&=&0,\\
\frac{\partial (v-v_0)}{\partial z}|_{z=0}=0,\quad \frac{\partial (v-v_0)}{\partial z}|_{z=-1}=0, \quad (v-v_0)\cdot n|_{\Gamma_l}=0, \quad \frac{\partial (v-v_0)}{\partial n}\times n|_{\Gamma_l}&=&0,\\
\nabla \Big(q_s-\int^z_{-1} \eta dz'\Big)+fk\times \gamma+L_1 \gamma&=&0,\\
\frac{\partial \gamma}{\partial z}|_{z=0}=0,\quad \frac{\partial \gamma}{\partial z}|_{z=-1}=0, \quad \gamma\cdot n|_{\Gamma_l}=0, \quad \frac{\partial \gamma}{\partial n}\times n|_{\Gamma_l}&=&0,
\end{eqnarray*}
by Proposition 2.5 in \cite{C-T-2}, we have
\begin{eqnarray}\label{eq-49}
|v-v_0|_{H^1(\mathcal{O})}&\leq& c|\chi-\chi_0|, \quad |v-v_0|_{H^2(\mathcal{O})}\leq c\|\chi-\chi_0\|,\\
\label{eq-50}
|\gamma|_{H^1(\mathcal{O})}&\leq& c|\eta|, \quad |\gamma|_{H^2(\mathcal{O})}\leq c\|\eta\|.
\end{eqnarray}
For $\eta$, we have
\begin{eqnarray}\label{eq-53}
\frac{d \eta}{dt}+L_2 \eta &-&(\delta\cdot \nabla)(\chi-\chi_0)+(v_0\cdot \nabla)\eta+((v-v_0)\cdot \nabla)(\chi-\chi_0)+(\gamma\cdot \nabla)\chi_0-((v-v_0)\cdot \nabla)\theta\\ \notag
&&-\Phi(\delta)\frac{\partial(\chi-\chi_0)}{\partial z}+\Phi(v_0)\frac{\partial \eta}{\partial z}+\Phi(v-v_0)\frac{\partial(\chi-\chi_0)}{\partial z}+\Phi(\gamma)\frac{\partial \chi_0}{\partial z}-\Phi(v-v_0)\frac{\partial \theta}{\partial z}=0.
\end{eqnarray}
Taking the inner product of equation (\ref{eq-53}) with $L_2 \eta$ in $H$, we have
\begin{eqnarray*}
\frac{1}{2}\frac{d \|\eta\|^2}{dt}+|L_2 \eta|^2
 &\leq&
 \Big|\int_{\mathcal{O}}(\delta\cdot \nabla)(\chi-\chi_0)L_2 \eta dxdydz\Big|+\Big|\int_{\mathcal{O}}\Phi(\delta)\frac{\partial(\chi-\chi_0)}{\partial z}L_2 \eta dxdydz\Big|
 \\
 &&+\Big|\int_{\mathcal{O}}(v_0\cdot \nabla)\eta L_2 \eta dxdydz\Big|+\Big|\int_{\mathcal{O}}\Phi(v_0)\frac{\partial \eta}{\partial z}L_2 \eta dxdydz\Big|
\\&&+\Big|\int_{\mathcal{O}}((v-v_0)\cdot \nabla)(\chi-\chi_0) L_2 \eta dxdydz\Big|  +\Big|\int_{\mathcal{O}}\Phi(v-v_0)\frac{\partial(\chi-\chi_0)}{\partial z}L_2 \eta dxdydz\Big|\\
 &&+\Big|\int_{\mathcal{O}}(\gamma\cdot \nabla)\chi_0 L_2 \eta dxdydz\Big|
+\Big|\int_{\mathcal{O}}\Phi(\gamma)\frac{\partial \chi_0}{\partial z}L_2 \eta dxdydz\Big|\\
&&+\Big|\int_{\mathcal{O}}((v-v_0)\cdot \nabla)\theta L_2 \eta dxdydz\Big|+\Big|\int_{\mathcal{O}}\Phi(v-v_0)\frac{\partial \theta}{\partial z}L_2 \eta dxdydz\Big|.
\end{eqnarray*}
We will estimate these terms one by one. By H\"{o}lder inequality and Lemma \ref{lem-2},
\begin{eqnarray*}
&&\Big|\int_{\mathcal{O}}(\delta\cdot \nabla)(\chi-\chi_0)L_2 \eta dxdydz\Big|+\Big|\int_{\mathcal{O}}\Phi(\delta)\frac{\partial(\chi-\chi_0)}{\partial z}L_2 \eta dxdydz\Big|\\
&\leq& c |L_2 \eta|\|\chi-\chi_0\|^{\frac{1}{2}}|L_2(\chi-\chi_0)|^{\frac{1}{2}}|\delta|_{H^1(\mathcal{O})}+c |L_2 \eta|\|\chi-\chi_0\|^{\frac{1}{2}}|L_2(\chi-\chi_0)|^{\frac{1}{2}}|\delta|^{\frac{1}{2}}_{H^1(\mathcal{O})}|\delta|^{\frac{1}{2}}_{H^2(\mathcal{O})},
\end{eqnarray*}
similarly, we have
\begin{eqnarray*}
&&\Big|\int_{\mathcal{O}}(v_0\cdot \nabla)\eta L_2 \eta dxdydz\Big|+\Big|\int_{\mathcal{O}}\Phi(v_0)\frac{\partial \eta}{\partial z}L_2 \eta dxdydz\Big|\\
&\leq& c|L_2 \eta|^{\frac{3}{2}}\|\eta\|^{\frac{1}{2}}|v_0|_6+c|L_2 \eta||\frac{\partial \eta}{\partial z}|^{\frac{1}{2}}|\nabla\frac{\partial \eta}{\partial z}|^{\frac{1}{2}}| v_0|^{\frac{1}{2}}_{H^1(\mathcal{O})}|v_0|^{\frac{1}{2}}_{H^2(\mathcal{O})},
\end{eqnarray*}
moreover, we get
\begin{eqnarray*}
&&\Big|\int_{\mathcal{O}}((v-v_0)\cdot \nabla)(\chi-\chi_0) L_2 \eta dxdydz\Big|  +\Big|\int_{\mathcal{O}}\Phi(v-v_0)\frac{\partial(\chi-\chi_0)}{\partial z}L_2 \eta dxdydz\Big|\\
&\leq& c|L_2 \eta|\|\chi-\chi_0\|^{\frac{1}{2}}|L_2(\chi-\chi_0)|^{\frac{1}{2}}|v-v_0|_6+c|L_2 \eta|\|\chi-\chi_0\|^{\frac{1}{2}}|L_2(\chi-\chi_0)|^{\frac{1}{2}}|v-v_0|^{\frac{1}{2}}_{H^1(\mathcal{O})}
|v-v_0|^{\frac{1}{2}}_{H^2(\mathcal{O})}\\
&\leq& c|L_2 \eta|\|\chi-\chi_0\|^{\frac{1}{2}}|L_2(\chi-\chi_0)|^{\frac{1}{2}}|\chi-\chi_0|+c|L_2 \eta|\|\chi-\chi_0\||L_2(\chi-\chi_0)|^{\frac{1}{2}}|\chi-\chi_0|^{\frac{1}{2}}
,
\end{eqnarray*}
where (\ref{eq-49}) is used.
\begin{eqnarray*}
&&\Big|\int_{\mathcal{O}}(\gamma\cdot \nabla)\chi_0 L_2 \eta dxdydz\Big|
+\Big|\int_{\mathcal{O}}\Phi(\gamma)\frac{\partial \chi_0}{\partial z}L_2 \eta dxdydz\Big|\\
&\leq& c|L_2 \eta||\nabla\chi_0|^{\frac{1}{2}}|L_2\chi_0|^{\frac{1}{2}}|\gamma|_6+c|L_2 \eta||\frac{\partial\chi_0}{\partial z}|^{\frac{1}{2}}|\nabla\frac{\partial\chi_0}{\partial z}|^{\frac{1}{2}}|L_2\chi_0|^{\frac{1}{2}}|\gamma|^{\frac{1}{2}}_{H^1(\mathcal{O})}|\gamma|^{\frac{1}{2}}_{H^2(\mathcal{O})}\\
&\leq& c|L_2 \eta||\nabla\chi_0|^{\frac{1}{2}}|L_2\chi_0|^{\frac{1}{2}}|\eta|+c|L_2 \eta||\frac{\partial\chi_0}{\partial z}|^{\frac{1}{2}}|L_2\chi_0||\eta|^{\frac{1}{2}}\|\eta\|^{\frac{1}{2}},
\end{eqnarray*}
where $|\gamma|_6\leq c\|\gamma\|$, and (\ref{eq-50}) is used.

At last, by (\ref{eq-49}), we obtain
\begin{eqnarray*}
&&\Big|\int_{\mathcal{O}}((v-v_0)\cdot \nabla)\theta L_2 \eta dxdydz\Big|+\Big|\int_{\mathcal{O}}\Phi(v-v_0)\frac{\partial \theta}{\partial z}L_2 \eta dxdydz\Big|\\
&\leq& c|L_2 \eta||\nabla\theta|^{\frac{1}{2}}|L_2 \theta|^{\frac{1}{2}}|v-v_0|_6+c|L_2 \eta||\frac{\partial\theta}{\partial z}|^{\frac{1}{2}}|\nabla\frac{\partial\theta}{\partial z}|^{\frac{1}{2}}|L_2 \theta|^{\frac{1}{2}}|v-v_0|^{\frac{1}{2}}_{H^1(\mathcal{O})}|v-v_0|^{\frac{1}{2}}_{H^2(\mathcal{O})}\\
&\leq& c|L_2 \eta||\nabla\theta|^{\frac{1}{2}}|L_2 \theta|^{\frac{1}{2}}|\chi-\chi_0|+c|L_2 \eta||\frac{\partial\theta}{\partial z}|^{\frac{1}{2}}|L_2 \theta||\chi-\chi_0|^{\frac{1}{2}}\|\chi-\chi_0\|^{\frac{1}{2}},
\end{eqnarray*}
Then, by the Young inequality, we get
\[
\frac{d\|\eta\|^2}{dt}+|L_2 \eta|^2\leq cg(t) \|\eta\|^2+ch(t),
\]
where
\[
g(t)=|v_0|^4_{H^1(\mathcal{O})}+|v_0|^2_{H^1(\mathcal{O})}|v_0|^2_{H^2(\mathcal{O})}
+|\chi_0|^2_{H^1(\mathcal{O})}|\chi_0|^2_{H^2(\mathcal{O})},
\]
and
\begin{eqnarray*}
h(t)&=&|\delta|^4_{H^1(\mathcal{O})}+|\delta|^2_{H^1(\mathcal{O})}|\delta|^2_{H^2(\mathcal{O})}+\|\theta\|^2
|L_2 \theta|^2+\|\chi-\chi_0\|^2|L_2(\chi-\chi_0)|^2\\
&&+|\chi-\chi_0|^2+|\chi-\chi_0|^2\|\chi-\chi_0\|^2
+\|\theta\|^2\|\chi-\chi_0\|^2.
\end{eqnarray*}
Thanks to (\ref{eq-51}), (\ref{eq-52}) and by Gronwall inequality, we have
\begin{eqnarray*}
\|\eta(t)\|^2
&\leq& \int^t_{t_0}h(s)ds+\int^t_{t_0}(\int^s_{t_0}h(\tau)d\tau)g(s) \exp (\int^t_s g(r)dr)ds\\
&\leq& \Big[c(1+t)+\int^t_{t_0}\|\chi-\chi_0\|^2|L_2(\chi-\chi_0)|^2ds+\int^t_{t_0}|\chi-\chi_0|^2ds\\
&&+\int^t_{t_0}|\chi-\chi_0|^2\|\chi-\chi_0\|^2ds
+\int^t_{t_0}\|\chi-\chi_0\|^2ds\Big]\exp (c+ct).
\end{eqnarray*}
Now, for $\chi-\chi_0$, using the similar method as above, we obtain
\[
\frac{d\|\chi-\chi_0\|^2}{dt}+|L_2(\chi-\chi_0)|^2\leq c\Big(\|\chi\|^2+\|\chi\||L_2\chi|+|v_0|^2_{H^1(\mathcal{O})}+|v_0|^2_{H^1(\mathcal{O})}|v_0|^2_{H^2(\mathcal{O})}\Big)\|\chi-\chi_0\|^2,
\]
by Theorem \ref{thm-1} and Gronwall inequality, it gives
\[
\|(\chi-\chi_0)(t)\|^2\leq C(t,\omega)\|\xi-\xi_0\|^2.
\]
Consequently, we have
\[
\int^t_{t_0}|L_2(\chi-\chi_0)|^2ds\leq C(t,\omega)\|\xi-\xi_0\|^2,
\]
thus,
\[
\|\eta(t)\|^2\leq C(t,\omega)\|\xi-\xi_0\|^4,
\]
where $C(t,\omega)$ is finite, $\mathbb{P}-$ a.s..
Recalling the definition of $g_{\delta}(t)$, we conclude it converges to zero in $L^{\infty}(\Omega, \mathbb{P})$ for each $t$. Collecting all the above facts, Theorem \ref{thm-8} is proved.

$\hfill\blacksquare$
\section{ Exponential Mixing }
 \subsection{Preliminaries on Coupling Method}
 Let $\lambda$ be a probability measure and $Y$ a random variable on a probability space $(\Omega,\mathcal{F},\mathbb{P})$. The law of $Y$ is denoted by $\mathcal{D}(Y)$.
 Given a Polish space $E$, the space $Lip_b(E)$ consists of all the bounded and Lipschitz real valued functions on $E$. Its norm is given by
 \[
 \|\phi\|_L=|\phi|_{\infty}+L_{\phi} \quad \phi\in Lip_b(E),
 \]
 where $|\cdot|_{\infty}$ is the sup norm and $L_{\phi}$ is the Lipschitz constant of $\phi$. The space of probability measures on $E$ is denoted by $\mathcal{P}(E)$. Recall that the total variation of a finite signed measure $\mu$ on $E$ is defined by
 \[
 \|\mu\|_{var}=\sup\Big\{|\mu(\Gamma)| | \Gamma\in \mathcal{B}(E)\Big\},
 \]
 where $\mathcal{B}(E)$ is the set of the Borelian subsets of $E$. It's well known that $\|\cdot\|_{var}$ is the dual norm of $|\cdot|_{\infty}$. The Wasserstein distance between probability measure $\mu_1$ and $\mu_2$ is defined by
 \begin{eqnarray}\label{eq-100}
 \|\mu_1-\mu_2\|_{*}=\sup_{\phi\in Lip_b(E),\|\phi\|_L\leq 1}\big|\int_{E}\phi(x)\mu_1(dx)-\int_{E}\phi(x)\mu_2(dx) \big|,
 \end{eqnarray}
 which is the dual norm of $\|\cdot\|_{L}$.

From Theorem \ref{thm-1}, we know that the solution $T$ is a Markov process on $H$ and depending measurably on the cylindrical wiener process $W$. We write
 \[
 T(t)=T(t,W,T_0),
 \]
 where $T_0=T(x,y,z,0)$ is the initial data $T(0,W,T_0)=T_0$.
Denote by $(\mathcal{P}_{t})_{t\in\mathbb{R}^{+}}$ the Markov transition semigroup associated to the Markov family $(T(\cdot,W,T_0))_{T_0\in H}$, which is defined by
\begin{equation}\label{eq-54}
\mathcal{P}_{t}f(x)=\mathbb{E}[f(T(t,x))],
\end{equation}
for $f\in B_b(H)$ and initial data $x\in H$.
 \subsection{A General Criterion for Exponential Mixing}
 In \cite{O-C}, the author established a general criterion which ensures exponential mixing of parabolic stochastic partial differential equations driven by multiplicative noise. In the following, we will describe the criterion and verify it holds for (\ref{eq-8})-(\ref{eq-14}).

 The basis of the criterion is the existence of an auxiliary process $\tilde{T} =\ \tilde{T}(t,W,T_0,\tilde{T_0})$, which verifying a set of conditions. The idea is that $\tilde{T}(t,W,T_0,\tilde{T_0})$ is close to $T(t,W,T_0)$ and has a law absolutely continuous with respect to the law of $T(t,W,\tilde{T_0})$. More precisely, Suppose that there exists a function
\[
\tilde{T}:\mathbb{R}^{+}\times C(\mathbb{R}^{+};\mathbb{R})^{\infty}\times H \times H\rightarrow H,
\]
satisfying the following assumptions.
\begin{description}
  \item[\textbf{A0}]  For each $T_0$, $\tilde{T}_0\in H$, $\tilde{T}(t,W,T_0,\tilde{T_0})$ is  non-anticipative and measurable map with respect to W. Moreover
      \[
       (T(t),\tilde{T}(t))=\Big(T(t,W,T_0), \tilde{T}(t,W,T_0,\tilde{T}_0)\Big)
       \]
defines an homogenous  Markov
process and its law $\mathcal{D}(T,\tilde{T})$  is measurably  with respect to $(T_0,\tilde{T_0})$.
\end{description}
The following assumptions involve a positive measurable functional $\mathcal{H}$ : $H\rightarrow$ $\mathbb{R}^{+}$, which plays the role of a Lyapunov functional.
\begin{description}
  \item[\textbf{A1}]  There exist  $c_{0}>0$,\ $C_{1}>0$ and for any $\alpha >0$, $C'_{\alpha}$  such that for any $T_0\in H$, any $t\geq0$, any $\alpha>0$ and any stopping time $\tau\geq 0$
 \begin{equation}
\left\{
\begin{array}{rcl}
\mathbb{E}\Big(\mathcal{H}(T(t,W,T_{0}))\Big) & \leq & e^{-c_{0}t}\mathcal{H}(T_{0})+C_{1},\\
\mathbb{E}\Big(e^{-\alpha \tau}\mathcal{H}(T(t,W,T_{0})I_{\tau<\infty}\Big)& \leq& \mathcal{H}(T_{0}) +C'_{\alpha}.%
\end{array}%
\right.
\end{equation}

  \item[\textbf{A2}] There exists $C>0$ such that for any $T^1_0$, $T^2_0\in H$ satisfying
  \begin{eqnarray}\label{equ-7}
  \mathcal{H}(T^1_0)+\mathcal{H}(T^2_0)\leq 2C_1,
  \end{eqnarray}
  then for any $W_1$,$W_2$ cylindrical Winner Processes on $H$ and for any $t\geq 0$, we have
  \begin{eqnarray*}
  &\mathbb{P}\Big(|T(t,W_2,T^2_0)-T(t,W_1,T^1_0)|\geq Ce^{-\gamma t}\ {\rm{and}} \  \tilde{T}(\cdot, W_1,T^1_0,T^2_0)=T(\cdot, W_2,T^2_0)\ {\rm{on}} \ [0,t]\Big)\\
  &\leq Ce^{-\gamma t}.
  \end{eqnarray*}

  \item[\textbf{A3}] There exists a function $h$: $H\times H\rightarrow H$ such that for any $(t,T^1_0,T^2_0)\in \mathbb{R}^{+}\times H\times H$ and $W$ cylindrical Winner process on $H$,  we have $\mathbb{P}$-a.s.
      \[
      \tilde{T}(t,W,T^1_0,T^2_0)=T\left(t,W+\int^{\cdot}_{0}h\big(T(s,W,T^1_0),\tilde{T}(s,W,T^1_0,T^2_0)\big)ds,T^2_0\right).
      \]
  \item[\textbf{A4}]  For any\ $W_1$,\ $W_2$ cylindrical Wiener processes on\ $H$, for any\ $T^1_0$,$T^2_0$ satisfying (\ref{equ-7}), for any $t_0\geq 0$, we have
\[
\mathbb{P}\left(\int^{\tau}_{t_0}|h(t)|^2dt\geq Ce^{-\gamma t_0}\ {\rm{and}}\  \tilde{T}(\cdot,W_1,T^{1}_{0},T^{2}_{0})=T(\cdot,W_2,T^{2}_{0})\ \rm{on}\ [0,\tau] \right)\leq Ce^{-\gamma t_0},
\]
\rm{where} $\tau\geq t_0$ \rm{is\ any\ stopping\ time\ and\ where}
\[
h(t)=h\Big(T(t,W,T^1_0),\tilde{T}(t,W,T^1_0,T^2_0)\Big).
\]
\item[\textbf{A5}] There exists\ $p_1>0$ such that for any\ $T^1_0$,\ $T^2_0$ in $H$ satisfying (\ref{equ-7}), we have
\[
\mathbb{P}\left(\int^{\infty}_{0}|h(t)|^2dt\leq C\right)\geq p_1,
\]
where
\[
h(t)=h\Big(T(t,W,T^1_0),\tilde{T}(t,W,T^1_0,T^2_0)\Big).
\]
\end{description}

\begin{thm}\cite{O-C}\label{theorem-1}
Under assumptions $\textbf{A0}-\textbf{A5}$, there exists a unique stationary probability measure $\mu$ of $(\mathcal{P}_t)_{t\in\mathbb{R}^{+}}$ on $H$. Moreover, $\mu$ satisfies
\[
\int_{H}\mathcal{H}(u)d\mu(u)< \infty,
\]
and there exist $C$, $\gamma'>0$ such that for any $\lambda \in \mathcal{P}(H)$
\[
\|\mathcal{P}^*_t\lambda-\mu\|_{*}\leq Ce^{-\gamma' t}\Big(1+\int_H \mathcal{H}(u) d \lambda(u)\Big).
\]
\end{thm}
\subsection{Main Proof of Theorem \ref{thm-3}}
We need to check the assumptions \textbf{A0-A5} for (\ref{eq-8})-(\ref{eq-14}). For this propers, we need  nondegenerate condition on the low modes. Such kind condition is necessary for the ergodicity of $\{{\cal{P}}_t\}_{t\in R^+},$ although it maybe not optimal in the following.
 \begin{description}
   \item[\textbf{Hypothesis H0}] There exists N$\in \mathbb{N}$ and a bounded measurable map $g\in \mathcal{L}(H,H)$ such that
   \[
   Gg=P_N.
   \]
 \end{description}
 \begin{remark}
 The above number $N$ will be chosen to be big enough in the following Proposition \ref{prp-4}.
 \end{remark}
Firstly, \textbf{A0} is the consequence of Theorem \ref{thm-1}.
Choosing $\mathcal{H}=|\cdot|^{2}$, the following Lemma \ref{lem-3} states that \textbf{A1} is true.
\begin{lemma}\label{lem-3}\textbf{(The Lyapunov Structure)} There exist  $c_{0}>0$, $C_{1}$ and a family $(C'_{\alpha})_{\alpha}$ only depending on $B_{0}$ and $D$, such that
\begin{equation}
\left\{
\begin{array}{rcl}
\mathbb{E}\Big(|T(t,W,T_{0})|^{2}\Big) & \leq & e^{-c_{0}t}|T_{0}|^{2}+C_{1},\\
\mathbb{E}\left(e^{-\alpha \tau}|T(t,W,T_{0})|^{2}I_{\tau<\infty}\right)& \leq& |T_{0}|^{2}+C'_{\alpha},%
\end{array}%
\right.
\end{equation}
 for any\ $t\geq0$, any\ $\alpha>0$ and any stopping time\ $\tau$.
 \end{lemma}
 \textbf{Proof of Lemma \ref{lem-3}} \quad
 Applying It\^{o}  formula to $|T|^{2}$,
 \begin{equation}\notag
 \frac{d|T|^{2}}{dt}+2\|T\|^2
\leq -2\left((v\cdot \nabla )T +\Phi(v)\frac{\partial T}{\partial z}, T\right)
+ 2(T,G\frac{dW}{dt})+B_0 ,
\end{equation}
where
\[
B_0=\|G\|^2_{\mathcal{L}_2(H,H)}.
\]
 By integration by parts, we have
\begin{equation}\notag
\int_{\mathcal{O}}[(v\cdot \nabla )T +\Phi(v)\frac{\partial T}{\partial z}]\cdot T dxdydz=0,
\end{equation}
then
\begin{equation}\notag
\frac{d|T|^{2}}{dt}+2\|T\|^2\leq 2(T,G\frac{dW}{dt})+B_0.
\end{equation}
Since $\|T\|^2\geq \lambda_1|T|^2$, applying It\^{o} formula to $e^{\lambda_1 t}|T|^2$, integrating and taking expectation, then we obtain
\begin{equation}\notag
\mathbb{E}|T(t)|^{2}\leq e^{-2\lambda_1 t}|T_0|^{2}+\frac{B_0}{\lambda_1},
\end{equation}
so we establish the first inequality of this lemma.

Let $\alpha >0$, applying It\^{o} formula to $e^{-\alpha t}|T|^{2}$ gives
\begin{eqnarray*}
 \frac{d(e^{-\alpha t}|T|^{2})}{dt}+2e^{-\alpha t}\|T\|^2
&\leq&-\alpha e^{-\alpha t}|T|^{2}-2e^{-\alpha t}\left(T,(v\cdot \nabla )T+\Phi(v)\frac{\partial T}{\partial z}-G\frac{dW}{dt}\right)+e^{-\alpha t}B_{0} \\
&\leq&-\alpha e^{-\alpha t}|T|^{2}+2e^{-\alpha t}\Big(T,G\frac{dW}{dt}\Big)+e^{-\alpha t}B_{0}.
\end{eqnarray*}
Let $\tau$ be a stopping time and $n\in \mathbb{N}$, similar to the above, it follows
by Gronwall Lemma,
\begin{equation}\notag
\mathbb{E}(e^{-\alpha \tau \wedge n }|T(\tau \wedge n)|^{2})\leq |T_{0}|^{2}+C'_{\alpha}.
\end{equation}
Let $n\rightarrow \infty $, we obtain the above inequality. $\hfill\blacksquare$

\

\textbf{
Construction of the auxiliary process}

Now, we construct the auxiliary process to make sure \textbf{A3} holds. Consider
\begin{equation}\label{eq-30}
\left\{
\begin{array}{rcl}
\nabla \tilde{p}_{s}-\int^{z}_{-1}\nabla \tilde{T}dz'+f\vec{k}\times \tilde{v}+L_1\tilde{v}
& = & 0,\\
\frac{\partial \tilde{T}}{\partial t}+(\tilde{v}\cdot\nabla)\tilde{T}+\Phi(\tilde{v})\frac{\partial \tilde{T}}{\partial z}+L_2\tilde{T}+P_N\delta(T,\tilde{T})
& = & G\frac{dW}{dt}.%
\end{array}%
\right.
\end{equation}
By \textbf{Hypothesis H0}, $Gg=P_N$, and set
\[
\tilde{W}_\cdot=W_\cdot+\int^{\cdot}_0-g(\tilde{T})\delta(T,\tilde{T})dt,
\]
then (\ref{eq-30}) can be written as
\begin{equation}\label{eq-31}
\left\{
\begin{array}{rcl}
\nabla \tilde{p}_{s}-\int^{z}_{-1}\nabla \tilde{T}dz'+f\vec{k}\times \tilde{v}+L_1\tilde{v}
& = & 0,\\
\frac{\partial \tilde{T}}{\partial t}+(\tilde{v}\cdot\nabla)\tilde{T}+\Phi(\tilde{v})\frac{\partial \tilde{T}}{\partial z}+L_2\tilde{T}
& = & G\frac{d\tilde{W}}{dt}.%
\end{array}%
\right.
\end{equation}
As we want $\tilde{T}$ and $T$ become very close in probability, it's natural to build $\tilde{T}$ such that (\ref{eq-30})  holds with
\begin{eqnarray}\notag
\delta(T,\tilde{T})=KP_N(\tilde{T}-T),
\end{eqnarray}
thus, \textbf{A3} holds with
\begin{equation}\label{equa-1}
h(T,\tilde{T})=-g(\tilde{T})\delta(T,\tilde{T})=-Kg(\tilde{T})P_N(\tilde{T}-T).
\end{equation}
In the following, we consider the following equation
\begin{equation}\notag
\left\{
\begin{array}{rcl}
\nabla \tilde{p}_{s}-\int^{z}_{-1}\nabla \tilde{T}dz'+f\vec{k}\times \tilde{v}+L_1\tilde{v}
& = & 0,\\
\frac{\partial \tilde{T}}{\partial t}+(\tilde{v}\cdot\nabla)\tilde{T}+\Phi(\tilde{v})\frac{\partial \tilde{T}}{\partial z}+L_2\tilde{T}+KP_N(\tilde{T}-T)
& = & G\frac{dW}{dt},%
\end{array}%
\right.
\end{equation}
where $K$ is a positive constant which will be chosen later and $N$ is the integer used in \textbf{Hypothesis H0}.

The following energy will play an important role in the remain part.
\begin{equation}\notag
\mathbb{E}_{T}(t)=|T(t)|^2+\int^{t}_{0}\|T(s)\|^2ds.
\end{equation}
\begin{lemma}\label{lem-2}
\textbf{(Exponential estimate  for the growth of the solution)} There exist positive constant $B_0$ and $\gamma_0=\gamma_0(B_0,\mathcal{O})$ such that
\[
\mathbb{E}\left(\exp\Big(\gamma_0\sup_{t\geq 0}(\mathbb{E}_{T}(t)-B_0 t)\Big)\right)\leq 2 e^{\gamma_0|T_0|^2},
\]
for any $T_0\in H$.
\end{lemma}
\textbf{Proof of Lemma \ref{lem-2}}  \quad
 We set for any $\gamma >0,$
 \[
 M(t)=2\int^{t}_{0}\Big(T(r),GdW(r)\Big), \quad \mathcal{M}_{\gamma}(t)=M(t)-\frac{\gamma}{2} \langle M\rangle (t).
 \]
 Remarking that
  \[
  d\langle M\rangle=4|G^{*} T|^2dt\leq 4cB_0\|T\|^{2}dt,
  \]
 where $c=\frac{1}{\lambda_1}$.
 From Lemma \ref{lem-3}, we have
 \begin{eqnarray}\notag
 d|T|^2+2\|T\|^2dt\leq dM(t)+B_0dt.
 \end{eqnarray}
 Setting $\gamma_1=\frac{1}{2cB_0}$, then
  \begin{equation}\label{eq-32}
  \mathbb{E}_{T}(t)\leq \mathcal{M}_{\gamma_1}(t)+|T_0|^2+B_0t.
  \end{equation}
Since $\gamma_1\mathcal{M}_{\gamma_1}(t)$ is a martingale, so $e^{\gamma_1\mathcal{M}_{\gamma_1}}$is a positive supermartingale whose value is 1 at time 0.
We deduce from maximal supermartingale inequality that
\[
\mathbb{P}\left(\sup_{t\geq0}\mathcal{M}_{\gamma_1}(t)\geq \rho\right)\leq \mathbb{P}\left(\sup_{t\geq0}e^{\gamma_1 \mathcal{M}_{\gamma_1}(t)}\geq e^{\gamma_1 \rho}\right)\leq e^{-\gamma_1 \rho}.
\]
Setting $\gamma_0=\frac{\gamma_1}{2}$. Notice that
 \[
 \mathbb{E}(e^{\gamma_0\sup\mathcal{M}_{\gamma_1}})=1+\gamma_0\int^{\infty}_{0}e^{\gamma_0 x}\mathbb{P}\Big(\sup\mathcal{M}_{\gamma_1}\geq x\Big)dx,
 \]
 which yields,
 \begin{eqnarray}\label{eq-33}
 \mathbb{E}(e^{\gamma_0\sup\mathcal{M}_{\gamma_1}})\leq 2.
 \end{eqnarray}
 Combining (\ref{eq-32}) and (\ref{eq-33}), it follows
\[
\mathbb{E}\left(\exp\Big(\gamma_0\sup_{t\geq 0}\big(\mathbb{E}_{T}(t)-B_0 t\big)\Big)\right)\leq 2 e^{\gamma_0|T_0|^2}.
\]
$\hfill\blacksquare$

Now, we are ready to prove \textbf{A2, A4} and \textbf{A5}.
\begin{prp}\label{prp-4}
 Assume \textbf{Hypothesis H0} holds. There exist $\gamma>0$,\ $\varepsilon\in(0,1]$,\ $K_0>0$ and\ $N_0=N_0(B_0,\mathcal{O})$ such that for any\ $K> K_0$ and $N>N_0$ and any $(t,T_0, \tilde{T_0})\in(0,\infty)\times H\times H$
\[
 \mathbb{E}\left(\Big(e^{t}|r(t)|^{2}+\frac{1}{2}\int^{t}_{0}e^{s}\| r(s)\|^{2}ds\Big)^{\varepsilon}\right)
 \leq 2|r(0)|^{2\varepsilon}e^{\gamma|T_0|^2},
\]
where $r(t)=\tilde{T}(t)-T(t)$.
Notice that, by Chebyshev inequality and (\ref{equa-1}), Proposition \ref{prp-4} obviously implies \textbf{A2, A4} and\textbf{ A5 }hold.
\end{prp}
\textbf{Proof of Proposition \ref{prp-4}} \quad  Let $u(t)=\tilde{v}(t)-v(t)$, and taking the difference between  (\ref{eq-8}), (\ref{eq-9}) and (\ref{eq-30}), we get
\begin{eqnarray*}
f{k}\times u+\nabla p_b-\nabla(\int^{z}_{-1} r dz')+L_1u&=&0,\\
\frac{\partial r}{\partial t}+\tilde{v}\cdot\nabla r+u\cdot \nabla T +\Phi(\tilde{v})\frac{\partial r}{\partial z}
+\Phi(u)\frac{\partial T}{\partial z}
+L_2 r+KP_N r& =&0,
\end{eqnarray*}
where $p_b=\tilde{p}_s-p_s$.

Applying It\^{o} formula to $|r|^2$, we get
\begin{equation}\notag
\begin{split}
 \frac{1}{2}\frac{d |r|^{2}}{d t}+\|r\|^2
 =&-\int_{\mathcal{O}}\big[\tilde{v}\cdot \nabla r+u\cdot \nabla T
 -(\int^{z}_{-1}\nabla \cdot \tilde{v}(x,y,z',t)dz')\frac{\partial r}{\partial z}\\
&\quad\quad-(\int^{z}_{-1}\nabla \cdot u(x,y,z',t)dz')\frac{\partial T}{\partial z}\big] \cdot r dxdydz
 -K|P_N r|^2,
\end{split}
\end{equation}
by integration by parts and boundary conditions, we get
\begin{equation}\notag
 \int_{\mathcal{O}}\big[\tilde{v}\cdot \nabla r-(\int^{z}_{-1}\nabla \cdot \tilde{v}(x,y,z',t)dz')\frac{\partial r}{\partial z}\big] \cdot r dxdydz=0,
\end{equation}
thus
\begin{eqnarray*}
\frac{1}{2}\frac{d |r|^{2}}{d t}+\|r\|^2+K|P_N r|^2 =-\int_{\mathcal{O}}(u\cdot \nabla) T r dxdydz
+\int_{\mathcal{O}}\Big(\int^{z}_{-1}\nabla \cdot u(x,y,z',t)dz'\Big)\frac{\partial T}{\partial z} r dxdydz.
\end{eqnarray*}
By similar method as Proposition \ref{prp-2}, we have
\begin{eqnarray}\notag
|u|_{H^1(\mathcal{O})}\leq c |r|,\quad |u|_{H^2(\mathcal{O})}\leq c \|r\|,
\end{eqnarray}
then, by Lemma \ref{lem-5},
\begin{eqnarray*}
\Big|\int_{\mathcal{O}} u\cdot \nabla T rdxdydz\Big|&\leq&\|T\||u|_{6}|r|_3\leq \|T\||u|_{H^1(\mathcal{O})}|r|_3 \leq c\|T\||r|^{\frac{3}{2}}\|r\|^{\frac{1}{2}},\\
\Big|\int_{\mathcal{O}}(\int^{z}_{-1}\nabla \cdot u(x,y,z',t)dz')\frac{\partial T}{\partial z} r dxdydz\Big| &\leq& c\|T\||u|^{\frac{1}{2}}_{H^1(\mathcal{O})}|u|^{\frac{1}{2}}_{H^2(\mathcal{O})}\|r\|^{\frac{1}{2}}|r|^{\frac{1}{2}}
\leq c\|T\|\|r\||r|,
\end{eqnarray*}
collecting all inequalities in the above and by the Young inequality, we have
\begin{equation}\label{equation2}
\frac{d|r|^2}{dt}+\frac{3}{2}\|r\|^2+K|P_N r|^2\leq c(1+\|T\|^2)|r|^2.
\end{equation}
Noticing that
\[
(K\wedge \mu_{N+1})|r|^2\leq |\nabla Q_N r|^2+|\partial_z Q_N r|^2+\alpha|Q_N r|_{z=0}|^2 +K|P_N r|^2\leq \|r\|^2+K|P_N r|^2,
\]
from (\ref{equation2}), we deduce
\[
d|r|^2+(K\wedge \mu_{N+1}-c-c\|T\|^2)|r|^2dt+\frac{1}{2}\|r\|^2dt\leq 0.
\]
Integrating this equality, it follows
\[
\mathbb{E}\left(e^t H(t)^{-1}|r(t)|^2+\frac{1}{2}\int^{t}_{0}e^sH(s)^{-1}\|r(s)\|^2ds\right)\leq |r(0)|^2,
\]
where
\[
H(t)=e^{-(K\wedge \mu_{N+1}-c-1)t+c\int^{t}_{0}\|T(s)\|^2ds}.
\]
From H\"{o}lder inequality, we get
\begin{eqnarray*}
&&\mathbb{E}\left(\Big(e^t|r(t)|^2+\frac{1}{2}\int^{t}_{0}e^s\|r(s)\|^2ds\Big)^{\varepsilon}\right)\\
&&\leq\sqrt{\mathbb{E}\Big(\sup_{(0,t)} H^{2\varepsilon}\Big)}
 \times  \left(\mathbb{E}\Big(e^t H(t)^{-1}|r(t)|^2+\frac{1}{2}\int^{t}_{0}e^sH(s)^{-1}\|r(s)\|^2ds\Big)\right)^{\varepsilon}.
\end{eqnarray*}
Choosing $N$ and $K$ large enough and $\varepsilon >0$ sufficiently small, it follows from Lemma \ref{lem-2} that
\[
\mathbb{E}\Big(\sup_{(0,t)} H^{2\varepsilon}\Big)\leq 2e^{\gamma_0|T_0|^2},
\]
which yields Proposition \ref{prp-4}. $\hfill\blacksquare$

\textbf{Proof of Theorem \ref{thm-3}} \quad From the above, we have established \textbf{A0-A5} hold.  This theorem is proved by Theorem \ref{theorem-1}.
$\hfill\blacksquare$

From (\ref{eq-100}), the definition of $\|\cdot\|_{*}$, we have
\begin{cor}\label{cor-1}
The Markov family $\{T(t,\omega)x\}$ is a system of mixing type in the following sense: it has a unique stationary measure $\mu$, such that
\[
\mathbb{E}f(T(t,\cdot)x)\rightarrow \int_H f(x)\mu(dx)\quad as \ t\rightarrow \infty,
\]
for each $f\in L(H)$ and each initial point $x\in H$.
\end{cor}

\section{Relationship Between Global Attractor and Invariant Measure}%
From Sects. 5 and 6, we have obtain the existence of global attractor $\mathcal{A}(\omega)$ and the unique invariant measure $\mu$ for the system (\ref{eq-8})-(\ref{eq-14}). In this section, we devote to exploring the relationship between $\mathcal{A}(\omega)$ and $\mu$.

We recall some definitions refer to \cite{RDS}.
\begin{dfn} Given an RDS $\varphi$. Then the mapping
\[
(\omega,x)\rightarrow (\theta(t)\omega, \varphi(t,\omega)x)=:\Theta(t)(\omega,x),\ t\in \mathbb{T},
\]
is a measurable DS on $(\Omega\times X, \mathcal{F}\otimes\mathcal{B}(X))$ which is called the skew product of the metric DS $(\Omega, \mathcal{F}, \mathbb{P},(\theta_t)_{t\in \mathbb{T}})$ and the cocycle $\varphi(t,\omega)$ on $X$.
\end{dfn}
\begin{dfn}
 Given a measurable RDS $\varphi$ over a metric DS $(\Omega, \mathcal{F},\mathbb{P},(\theta_t)_{t\in \mathbb{T}})$, a probability $\tilde{\mu}$ on $(\Omega\times X, \mathcal{F}\otimes\mathcal{B}(X) )$ is said to be an invariant measure for $\varphi$, if it satisfies
\begin{itemize}
  \item $\Theta(t) \tilde{\mu}=\tilde{\mu}$,
  \item $\pi_{\Omega}\tilde{\mu}={\mathbb{P}}$.
\end{itemize}
\end{dfn}
Denote $\mathcal{P}_{\mathbb{P}}$ = $\Big\{  \tilde{\mu} $ probability measure on $ (\Omega \times X,\mathcal{F}\otimes\mathcal{B}(X))$ with marginal $\mathbb{P}$ on $(\Omega,\mathcal{F})\Big\}$.
\begin{dfn}
\textbf{(Factorization of $\tilde{\mu}$)} Suppose $\tilde{\mu}\in \mathcal{P}_{\mathbb{P}}(\Omega \times X)$. We call a function $\mu_{\cdot}(\cdot): \Omega \times \mathcal{B}(X)\rightarrow[0,1]$ a factorization of $\tilde{\mu}$ with respect to $\mathbb{P}$, if
\begin{enumerate}
  \item for all $B\in \mathcal{B}, \omega\rightarrow \mu_{\omega}(B)$ is $\mathcal{F}$-measurable.
  \item for $\mathbb{P}$-a.s. $B\rightarrow \mu_{\omega}(B)$ is a probability measure on $(X,\mathcal{B}(X))$.
  \item for all $ A \in \Omega \times \mathcal{B}(X)$
  \[
  \tilde{\mu}(A)=\int_{\Omega}\int_{X}1_A(\omega,x)\mu_{\omega}(dx) \mathbb{P}(d\omega).
  \]
\end{enumerate}
\end{dfn}
Refer to Proposition 1.4.3 in \cite{RDS}, for RDS $\varphi$ over a metric DS $(\Omega, \mathcal{F},\mathbb{P},(\theta_t)_{t\in \mathbb{T}})$ on a Polish space X, $\tilde{\mu}$ is an invariant measure of $\varphi$, then a factorization of $\tilde{\mu}$ exists and is unique, $\mathbb{P}$-a.e.  We write symbolically
\[
\tilde{\mu}(d\omega, dx)={\mu}_{\omega}(dx)\mathbb{P}(d\omega).
\]


From Theorem \ref{thm-3}, we know there is a unique invariant measure $\mu$ for  the semigroup $\{\mathcal{P}_t\}$ associated with (\ref{eq-8})-(\ref{eq-14}), it implies that there exists a unique invariant measure for the RDS generated by (\ref{eq-8})-(\ref{eq-14}).
In fact, refer to \cite{SA}, we have
\begin{prp}\label{prp-6}
There exits a unique invariant measure $\tilde{\mu}$ for the skew product $\Theta(t)$, which is given by
\begin{eqnarray}\label{eq-55}
\mu_{\omega}=\lim_{t\rightarrow \infty}\varphi(t,\theta_{-t}\omega)\mu,
\quad \tilde{\mu}(d\omega, dx)={\mu}_{\omega}(dx)\mathbb{P}(d\omega).
\end{eqnarray}
\end{prp}
%

We now discuss the relationship between invariant measures and random attractors for the system (\ref{eq-8})-(\ref{eq-14}). We denote $\mu_{\omega}$ is the disintegration of the above $\tilde{\mu}$ and set
\begin{eqnarray*}
\mathcal{A}(\omega)=\left\{
                      \begin{array}{ll}
                        supp \mu_{\omega}, & \omega \in \Omega_0, \\
                        H, & \omega \notin \Omega_0,
                      \end{array}
                    \right.
\end{eqnarray*}
where $\Omega_0 \in \mathcal{F}$ is a set of full measure on which the limit in (\ref{eq-55}) exists.
\begin{cor}\label{thm-9}
Suppose $\tau\in L^2(D)$ and \textbf{Hypothesis H0} holds, then $\mathcal{A}(\omega)$ defined above is a minimal random attractor for (\ref{eq-8})-(\ref{eq-14}). Moreover, $supp \mu_{\omega}$ has finite Hausdorff dimension, $\mathbb{P}-a.s.$.
\end{cor}
\textbf{Proof of Corollary\ref{thm-9}}\quad  
To prove $\mathcal{A}(\omega)$ defined above is a minimal random attractor for (\ref{eq-8})-(\ref{eq-14}), refer to \cite{SA}, we only need to verify two properties: mixing property and compactness.
It's easy to know Corollary \ref{cor-1} is exactly the required mixing property and Theorem \ref{thm-2} gives the compactness property. Furthermore, by Theorem \ref{thm-8} and the general result $supp  \mu_{\omega} \subset\mathcal{A}(\omega) $ a.s, we conclude that $supp \mu_{\omega}$ has finite Hausdorff dimension.
$\hfill\blacksquare$

\vskip0.5cm {\small {\bf  Acknowledgements}\ \  This work was supported by National Natural Science Foundation of China (NSFC) (No. 11431014, No. 11371041, No. 11271356)},
Key Laboratory of Random Complex Structures and Data Science, Academy of Mathematics and Systems Science, Chinese Academy of Sciences(No. 2008DP173182).


\def\refname{ References}

\end{document}